\def\Date{July 28, 2011}

\magnification=\magstephalf
\hsize=14.5true cm
\vsize=24.6true cm
\frenchspacing\parskip4pt plus 1pt
\lineskiplimit=-2pt

\newread\aux\immediate\openin\aux=\jobname.aux 
\ifeof\aux\message{ <<< Run TeX a second time >>> }
\batchmode\else\input\jobname.aux\fi\closein\aux

\input nil.def\immediate\openout\aux=\jobname.aux

\openup-.18pt

\draftfalse 

\ifdraft\footline={\hss\sevenrm\Date\hss}\else\nopagenumbers\fi

\headline={\ifnum\pageno>1\eightrm\ifodd\pageno\hfill~~ and affinely
homogeneous surfaces \hfill{\tenbf\folio}\else{\tenbf\folio}\hfill
Nilpotent algebras \hfill\fi\else\hss\fi}

\centerline{\Gross Nilpotent algebras}
\medskip \centerline{\Gross and affinely homogeneous surfaces}

\bigskip\bigskip \centerline{By {\sl Gregor Fels} and {\sl Wilhelm
Kaup}} {\parindent0pt\footnote{}{\ninerm 2000 Mathematics Subject
Classification: 14J70, 32V40}}

\KAP{Introduction}{Introduction}

This paper is devoted to the investigation of finite dimensional
commutative nilpotent (associative) algebras $\5N$ over an arbitrary
base field $\FF$ of characteristic zero.  Our main attention is
focused on those algebras which have 1-dimensional annihilator since
these algebras naturally occur in connection with various,
geometrically motivated problems.  The unital extensions
$\5N^{\,0}=\FF\oplus\5N$ of such algebras are exactly the Gorenstein
algebras of finite positive vector space dimension over $\FF$.

There is only very sparse literature concerning the structure of
general nilpotent commutative algebras. For example, every such
algebra has a realization as a subalgebra of some $\End(V)$, $V$ a
vector space over $\FF$, which is maximal with respect to the property
that it consists of nilpotent and commuting endomorphisms.  This means
that abstractly given nilpotent algebras and nilpotent subalgebras of
endomorphism algebras are essentially the same. While the structure
and classification of maximal commutative algebras consisting of
semisimple endomorphisms (Cartan subalgebras) is very well understood,
there is only little known about the general structure, not to mention
a classification, of its nilpotent counterpart (for an approach in
terms of Macaulay's inverse systems compare with \Lit{ELRO}). This is a
bit surprising, as such nilpotent algebras are quite ubiquitous
objects which occur in various areas of mathematics.  One reason is
certainly that the theory of nilpotent algebras is more involved than
the theory of Cartan subalgebras, due to the lack of rigidity
properties and other obvious visible structure.  Since the standard
tools from the Cartan theory such as the root theory cannot be applied
in the nilpotent case, the desire arises for appropriate objects which
help to understand nilpotent commutative algebras. One of the purposes
of this paper is to develop such tools.

From a purely algebraic point of view nilpotent commutative algebras
are building blocks for general commutative algebras which, for
instance, seems to be very important for quantum physics.  As already
mentioned, commutative nilpotent algebras naturally arise in the
context of several geometrically motivated questions as they often
serve as invariants attached to certain geometric objects. Our
interest in commutative nilpotent algebras also originates from
geometry.  To be more specific, we mention two types of geometric
problems which provide us with commutative nilpotent algebras which in
turn encode some of the geometric structure of the original questions.

In Cauchy-Riemann geometry there is the question under which
conditions certain CR-manifolds are (locally) equivalent to tube
manifolds $S\times i\RR^n\subset \RR^n\oplus i\RR^n=\CC^n$ and how
many different tube realizations do exist. In \Lit{FEKA} we show that
this geometric problem (in the case of non-degenerate hyperquadrics)
can be reduced to the classification of real {\sl and} complex
commutative nilpotent subalgebras with 1-dimensional annihilator.

Another type of problems arises from the study of isolated
hypersurface singularities and their versal deformations: Let $h$ be
(the germ of) a holomorphic function, defined in a neighbourhood $U$
of $0\in \CC^n$, i.e.$\,$, $h\in \5O_n:=\CC\{z_1,...\,, z_n\}$ such
that the hypersurface $\{h=0\}$ has an isolated singularity in
$0$. This implies that $\grad h:U\to \CC^n$ is a finite map, and
consequently
$$\qu{\5O_n}{h^*(\5O_n)\cd \5O_n}=\qu{\5O_n}{J(h)}$$ is a finite
dimensional local algebra. Here, $J(h):=\big({\partial h\over \partial
z_1},...\,, {\partial h\over\partial z_n}\big)$ denotes the {\sl
Jacobi ideal} of $h$ in $\5O_n$. As a consequence of Nakayama's Lemma
its maximal ideal is a (commutative) nilpotent algebra.  The full
algebra serves as parameter space for the universal deformation of
the isolated singularity of $h$. Consequently also the maximal ideal of
the Tjurina algebra $\qu{\5O_n}{(h,J(h))}$ is nilpotent and
its algebra structure turns out to determine the original isolated
singularity up to a biholomorphic equivalence, see e.g. \Lit{XUYA} for
further details.

For short, we call from now on a finite dimensional commutative
nilpotent algebra over $\FF$ with 1-dimensional annihilator simply an
{\sl admissible algebra}.  In this paper we give a construction of
objects naturally associated with an admissible algebra $\5N$, encoding
sufficient information to recover the original algebra. These objects
seem to be easier to deal with than the admissible algebras themselves
and also serve as convenient invariants allowing the explicit
verification of whether two admissible algebras are isomorphic.  These
objects are certain classes of smooth subvarieties of $\5N\cong
\FF^{n+1}$ as well as certain classes of polynomials, which we call
{\sl nil-polynomials}. These polynomials are closely related to the
aforementioned smooth subvarieties.  Roughly speaking the
nil-polynomials are certain truncated exponential series (here the
nilpotency is crucial), concatenated with a linear functional, which
essentially is nothing but a linear projection onto the annihilator of
$\5N$. The constant and linear parts of every nil-polynomial
$\6p\in\FF[X_{1},\dots,X_{n}]$ vanish, but the quadratic part of $\6p$
is non-degenerate.  Up to isomorphism the algebra structure on $\5N$
can be recovered from the polynomial $\6p$. Even more is true, as the
quadratic plus cubic term alone suffice to determine the structure of
$\5N$, and in turn the entire nil-polynomial $\6p$. Unless $n=0$, its
degree coincides with the nil-index of $\5N$, i.e., the maximal number
$\nu$ with $\5N^{\,\nu}\ne 0$.

Let $\5A$ be the annihilator of $\5N$ and $\5K$ a hyperplane in $\5N$
transversal to $\5A$, that is, $\5N=\5K\oplus\5A$.  The smooth variety
associated with $\5N$ (and depending on linear isomorphisms
$\5K\cong\FF^{n}$, $\5A\cong\FF$) is simply the graph
$S\subset\FF^{n+1}$ of the corresponding nil-polynomial $\6p:\FF^n\to
\FF$.  We call every such $S$ a {\sl nil-hypersurface.}  Among other things
we prove that two admissible algebras $\5N$, $\tilde\5N$ with
nil-hypersurfaces $S$, $\tilde S$ are isomorphic as algebras if and
only if $S$, $\tilde S$ are affinely equivalent. For an even stronger
statement see Theorem \ruf{AB}. We also show that affine equivalence
for $S$, $\tilde S$ can be replaced by linear equivalence if and only
if the nil-hypersurface $S$ is affinely homogeneous. Linear
equivalence gives a stronger and computationally more convenient
condition for the isomorphy of the algebras $\5N$, $\tilde\5N$. On the
polynomial level it means that for the corresponding nil-polynomials
$\6p$, $\tilde\6p$ there is a $g\in\GL(n,\FF)$ such that $\tilde\6p$
and $\6p\circ g$ differ by a constant factor from $\FF^{*}$.  We
further establish a duality between a fixed nil-hypersurface $S$ of
$\5N$ and the parameter space $\Sigma(\5N)$ of all such
nil-hypersurfaces.  Taking this duality a step further, we show that
the action of the affine group $\Aff(S)$ on $S$ is equivariantly
isomorphic to the action of the algebra automorphism group $\Aut(\5N)$
on the affine space $\Pi(\5N)$ of all projections with range the
annihilator of $\5N$.

As already mentioned, affine homogeneity of an associated
nil-hypersurface $S$ of $\5N$ makes computations more efficient.
However, the question for which admissible algebras $\5N$ the
nil-hypersurface $S$ is affinely homogeneous is quite involved.  Only
recently we were able to give a satisfactory answer to this question:
While the nil-hypersurface of every admissible algebra of nil-index
smaller than 5 is automatically affinely homogeneous, there are
non-homogenous counterexamples starting with nil-index 5.  In the case
however, where $\5N$ admits a $\ZZ^+$-gradation, every corresponding
nil-hypersurface $S$ is affinely homogeneous.

\medskip Our paper is organized as follows: In Section
\ruf{Preliminaries} we fix notation and state some preliminaries.  A
simple example is given, which indicates why in the rest of the paper
we stick to nilpotent algebras that are {\sl commutative}. In Section
\ruf{Affinely} we introduce the notion of a nil-surface, which is a
smooth algebraic variety $S_\pi$, associated with a given nilpotent
commutative algebra $\5N$, however depending also on a projection
$\pi\in \End(\5N)$.  Further we discuss various notions of gradations
for nilpotent commutative algebras.  The main result of the section
holds for $\5N$ admitting certain types of generalized gradations and
compatible projections $\pi$: In this case $S_{\pi}$ is affinely
homogeneous. A special version of this result for base field
$\FF\in\{\RR,\CC\}$ is already contained in \Lit{FEKA} and later also
was used in \Lit{FIKK} for $\FF=\CC$, compare also with \Lit{ISAE}. In
Section \ruf{Admissible} we restrict our attention to {\sl admissible
algebras} algebras $\5N$.  From that point on we only consider {\sl
admissible projections} $\pi$ on $\5N$ which means that the range of
$\pi$ is the (1-dimensional) annihilator of $\5N$. The main result of
this section is, roughly speaking, that the algebra structure of $\5N$
only depends on the hypersurface $S_\pi\subset \5N$, and `essential'
properties of $S_{\pi}$ do not depend on the admissible projection
$\pi$. Another statement is that the affine group $\Aff(S_{\pi})$ is
canonically isomorphic to the algebra automorphism group
$\Aut(\5N)$. The main result of the section, Theorem \ruf{AB}, is a
generalization and extension of a result in \Lit{FIKK} from base field
$\CC$ to arbitrary $\FF$. The proof in \Lit{FIKK} is of analytic
nature and we get the extension by applying a Lefschetz principle type
argument. We also investigate functorial properties of the space
$\Sigma(\5N)$ of all nil-hypersurfaces and of the affine space
$\Pi(\5N)$ of all admissible projections, formulate a duality
statement between a member $S$ of $\Sigma(\5N)$ and the family
$\Sigma(\5N)$ itself and prove the equivariant equivalence of the
natural actions of $\Aut(\5N)$ on $\Pi(\5N)$ and of $\Aff(S)$ on $S$.
We close the section with an infinitesimal analogon, more precisely,
we show for every base field of characteristic 0 that for any
admissible algebra $\5N$ with associated nil-hypersurface
$S\subset\5N$ the derivation algebra $\der(\5N)$ is isomorphic to the
Lie algebra $\aff(S)$ of all affine transformations $\5N\to\5N$ that
are `tangent' to $S$. In Section \ruf{Nil} we associate to every
admissible algebra $\5N\cong\FF^{n+1}$ a class of mutually affinely
equivalent polynomials in $\FF[x_{1},\dots,x_{n}]$, called {\sl
nil-polynomials.}  The graphs of these polynomials are affinely
equivalent to $S_{\pi}$.  In Section \ruf{Representations} we present
for every admissible $\5N$ certain canonical decompositions and show:
Every $\5N$ with nil-index $\le3$ has a grading, and for every $\5N$
with nil-index $\le4$ every $S_{\pi}$ is affinely homogeneous. Both
bounds for the nil-index are sharp as will be shown by counterexamples
in the last section. We also show for every admissible algebra $\5N$
that $\der(\5N)$ and $\Aut(\5N)$ have at least dimension
$\dim(\5N/\5N^{4})$.  In Section \ruf{examples} we give large classes
of admissible algebras of nil-index 3 and 4. It turns out in
particular, that in every dimension $\ge7$ for $\FF\in\{\RR,\CC\}$ the
number of isomorphy classes of admissible algebras is uncountable
infinite. We also get a classification of all admissible algebras of
nil-index 3. For the special case of algebraically closed base fields
this has already been achieved in \Lit{ELRO} by completely different
methods.  In Section \ruf{counter} we present various counterexamples
elaborated with computer aid. Among these we give an admissible
algebra $\5N$ of dimension 23 and nil-index 5 such that $S_{\pi}$ is
{\sl not affinely homogeneous.}  This disproves the Conjecture at the
end of \Lit{ISAA}, repeated and extended as Conjecture 2.4 in
[\ruf{ISAE}a].

One of the essential parts of the present paper is Theorem \ruf{NB}.
For the special case of {\sl admissible} algebras it also
occurs as Cor. 2.6 in [\ruf{ISAE}a] together with the statement {\sl
``We note that Corollary 2.6 was obtained by W. Kaup approximately
three months before this paper was written''} on page 3. This
Corollary is essentially the same as Theorem 2.5 in [\ruf{ISAE}a]. In
later versions such as [\ruf{ISAE}b] any hint to our priority is
missing.

\KAP{Preliminaries}{Preliminaries}

Throughout the paper $\NN$ is the set of all {\sl non-negative}
integers while $\ZZ^{+}$ is the semigroup of all {\sl positive}
integers. Further, $\FF$ is an arbitrary but fixed field of
characteristic 0. All algebras in the following are defined over $\FF$
and are associative, but may have infinite dimension as $\FF$-vector
spaces (at least in the first three sections). For every such algebra
$A$, every $x\in A$ and every integer $k\ge1$ we put
$$x^{(k)}:={1\over
k!}x^{k}\Steil{and}x^{(0)}:=\One\steil{if}A\steil{has a
unit}\One\,.\Leqno{KP}$$ Also we denote for every $j\in\ZZ^{+}$ by
$$\exp_{j}=\sum_{k=j}^{\infty}T^{(k)}\;\in\;\FF[[T]]$$ the $j$-{\sl
truncated exponential series.} Then
$\exp_{1}\!\circ\log_{1}=\log_{1}\!\circ\exp_{1}=T\,$ for
$$\log_{1}:=\sum_{k=1}^{\infty}{(-1)^{k+1}\over k}\,T^{k}\,.$$

\Definition{NA} For every $\FF$-algebra $A$ define inductively the
characteristic ideals $A^{k}\subset A$ by $A^{1}=A$ and
$A^{k+1}=\langle AA^{k}\rangle_{\FF}$. Also put $A_{[0]}=0$ and
$A_{[k]}:=\{x\in A:xA^{k}+A^{k}x=0\}$ for $k>0$.  Then $A$ is called
{\sl nilpotent} if $A^{k+1}=0$ for some $k\ge0$, and the minimal $k$
with this property is called the {\sl nil-index} of $A$, that we
denote by $\nu=\nu(A)$.

For nilpotent $A$ with nil-index $\nu$ the inclusion $A^{k}\subset
A_{[\nu+1-k]}$ is obvious for all $k\le\nu$ as well as
$\nu=\inf\{k\ge0:A_{[k]}=A\}$ . The ideal $\Ann(A):=A_{[1]}$, called
the {\sl annihilator} of $A$, plays a prominent role in the
following. The annihilator coincides with the {\sl socle} of $\5N$,
that is, the sum of all minimal ideals.

Let $\5N$ be a nilpotent algebra. From now on we always consider
$\exp_{1},\log_{1}:\5N\to\5N$ as polynomial mappings that are inverse
to each other. Fix an arbitrary projection
$\pi=\pi^{2}\in\End(\5N)$. Then for the polynomial mapping
$$\6f:\5N\to\5N\,,\qquad\6f(x):=\pi(\exp_{1}x)\;,\Leqno{DN}$$
$$S=S_{\pi}:=\6f^{-1}(0)=\log_{1}(\ker\pi)\;\subset\;\5N\Leqno{ND}$$
is a smooth algebraic subvariety of codimension $\rank(\pi)$
containing the origin. We will be mainly interested in the case where
the annihilator of $\5N$ has dimension 1 and is the range $\pi(\5N)$
of the projection $\pi$. Then $S_{\pi}$ is a hypersurface in $\5N$.

\smallskip With $\,\Aff(\5N)\,\cong\,\GL(\5N)\,\ltimes\,\5N$ we denote
the group of all affine bijections of $\5N$ and by
$\Aff(S):=\{g\in\Aff(\5N):g(S)=S\}$ the subgroup stabilizing
$S$. Furthermore, $\GL(S):=\{g\in\GL(\5N):g(S)=S\}$ is the isotropy
subgroup of $\Aff(S)$ at the origin. We are interested in cases where
$S$ is affinely homogeneous, that is, the group $\Aff(S)$ acts
transitively on $S$. This is not always true: As a counterexample
consider the matrix algebra $\5T$ of all strictly upper triangular
$n\times n$-matrices with coordinates $x_{jk}$ for $1\le j<k\le n$ and
projection $\pi$ given by $x\mapsto x_{1n}$ (after identifying
$\Ann(\5T)\cong\FF$ in the obvious way). For instance, for $n=4$ the
corresponding polynomial $\6f$ is given by
$$\6f(x)={1\over6}x_{12}x_{23}x_{34}+
{1\over2}(x_{13}x_{34}+x_{12}x_{24})+x_{14}$$ and it can be seen
that $S_{\pi}$ is not affinely homogeneous. Notice that the quadratic
part of $\6f$ is a {\sl degenerate} quadratic form on $\ker(\pi)$, in
contrast to the commutative case below.

 In case $\5N$ is commutative, for every projection $\pi$ on $\5N$
with range $\5N_{[1]}=\Ann(\5N)$ it is well known that for the
$\5N_{[1]}$-valued symmetric 2-form
$$\6b_{\pi}:\5N\times\5N\to\5N_{[1]},\qquad(x,y)\mapsto\pi(xy)\Leqno{NF}$$
the radical $\{x\in\5N:\6b_{\pi}(x,\5N)=0\}$ coincides with
$\5N_{[1]}$ (for a simple proof compare e.g. Prop. 2.1 in
\Lit{FIKK}). In particular, the form $\6b_{\pi}$ is non-degenerate on
$\pi^{-1}(0)\cong\5N/\5N_{[1]}$. The $\5N_{[1]}$-valued polynomial
$\6f=\pi\circ\exp_{1}$ has a unique decomposition
$\6f=\sum_{k\ge1}\6f^{[k]}$ into homogeneous components $\6f^{[k]}$ of
degree $k$. Clearly, $\6f^{[1]}=\pi$,
$\6f^{[2]}(x)={1\over2}\6b_{\pi}(x,x)$ and $\6f^{[\nu]}(x)=x^{(\nu)}$
for all $x\in\5N$ and $\nu=\nu(\5N)$.

\bigskip From now on we assume that $\5N$ is a {\sl commutative
nilpotent} algebra. In this paper we investigate properties of $\5N$,
the polynomials $\6f$ and the corresponding nil-surfaces in $\5N$ in a
fairly general algebraic setting. Our motivation however comes from
complex geometry.  For instance -- as already mentioned in the
introduction -- every real hyperquadric $Y$ in a complex projective
space $\PP_{n}(\CC)$ gives rise to several (real and complex)
commutative nilpotent algebras $\5N$.  Roughly speaking, the varieties
$S_\pi$ in case $\pi(\5N)=\5N_{[1]}$ occur as building blocks of bases
$F\subset \RR^n$ in various tube representations $F\times
i\RR^n\subset \CC^n$ of the CR-manifold $Y$, compare \Lit{FEKA}.
Another source of commutative nilpotent algebras arises in the context
of isolated hypersurface singularities, compare \Lit{FIKK}.

\medskip With $\5N^{\,0}:=\FF\cd\One\oplus\5N$ we denote the {\sl
unital extension} of $\5N$, having $\One$ as unit. This notation has
been chosen since then the canonical filtration extends as
$$\5N^{\,0}\supset\5N^{\,1}\supset\cdots\5N^{\,k}\supset\cdots
\Steil{with}\5N^{\,j}\5N^{\,k}\subset\5N^{\,j+k}\;,$$ $\5N^{1}=\5N$
and $\qu{\5N^{\,0}}{\5N^{\,1}}\cong\FF$.  Clearly, $\5N$ is the unique
maximal ideal of $\5N^{\,0}$. In case
$\5N=\5N_{1}\oplus\5N_{2}\oplus\dots\oplus\5N_{d}$ for linear
subspaces $\5N_{k}$ and a fixed integer $d\ge1$ we can write every
$x\in\5N$ as tuple $x=(x_{1},x_{2},\dots,x_{d})$ with
$x_{k}\in\5N_{k}$ and then
$$\exp(x)=\sum_{\mu\in\NN^{d}}x^{(\mu)}\;\in\;\One+\5N\;\subset\;\5N^{\,0},$$
where $x^{(\mu)}:=x_{1}^{(\mu_{1})}x_{2}^{(\mu_{2})}\cdots
x_{d}^{(\mu_{d})}$ for all $\mu=(\mu_{1},\dots,\mu_{d})$ and notation
\Ruf{KP} is in force.

For every nilpotent algebra $\5N$ the {\sl Hilbert
function} $H=H_{\5N}:\NN\to\NN$ of $\5N^{\,0}$ is defined by
$H(k)=\dim(\5N^{k}/\5N^{k+1})$. Clearly, $H(0)=1\le H(\nu)$ for
$\nu:=\nu(\5N)$, and $H(k)=0$ for $k>\nu$. As usual, we write $H$ also
as finite sequence $\{H(0),H(1),\dots,H(\nu)\}$ and call $H$ {\sl
symmetric} if $$H(k)=H(\nu-k)\steil{for all}0\le
k\le\nu\;.\Leqno{UH}$$ The Hilbert function is a rough invariant for
nilpotent algebras that in general does not suffice to distinguish two
given algebras.

The complement of $\5N$ in $\5N^{\,0}$ is the maximal subgroup of
$\5N^{\,0}$. A special subgroup is the unipotent group
$\,\5U:=\One+\5N$.  The exponential mapping $\exp:\5N\to\5U$ is a
group isomorphism with inverse $\log\,$, where
$\log(\One+x)=\log_{1}(x)$ for all $x\in\5N$.

With $\Aut(\5N)$ we denote the algebra automorphism group of
$\5N$. The endomorphism algebra $\End(\5N)$ endowed with the
commutator product $[\lambda,\sigma]=\lambda\,\sigma-\sigma\lambda$ is
a Lie algebra that we also denote by $\gl(\5N)$. With
$\der(\5N)\subset\gl(\5N)$ we denote the Lie subalgebra of all
derivations. This is an $\5N^{\,0}$-leftmodule in an obvious way. For
every nilpotent $\lambda\in\der(\5N)$ the operator
$\exp(\lambda)\in\Aut(\5N)$ is unipotent and conversely, every
unipotent $u\in\Aut(\5N)$ is of this form.

\medskip
Derivations of $\5N$ can be obtained in the following way: Let $\pi$
be a projection on $\5N$ with range $\5N_{[1]}$ and suppose that
$\lambda\in\End(\5N)$ satisfies $\lambda(\5N)\subset\5N_{[1]}$ and
$\lambda(\5N^{\,2})=0$. Then for every $a\in\5N_{[3]}$ the operator
$x\mapsto ax+\pi(ax)+\lambda(x)$ is in $\der(\5N)$. This implies
$$\dim\der(\5N)\;\ge\;(\dim\5N/\5N^{\,2})\dim\5N_{[1]}\,+\,
\dim\,(\5N_{[3]}/\5N_{[2]})\;.\Leqno{TY}$$ This estimate improves
Theorem 5.4 in \Lit{PERE}, which gives in case of algebraically closed
base field (and finite dimension) the lower bound \Ruf{TY} without the
summand $\dim\,(\5N_{[3]}/\5N_{[2]})$. Also
$$\dim\der(\5N)\;\le\;(\dim\5N/\5N^{\,2})\dim\5N$$ always holds since
every derivation of $\5N$ is uniquely determined by its values on a
linear subspace $L$ with $\5N=L+\5N^{2}$. The latter inequality is an
equality e.g. if $\5N=\7m\!/\!\7m^{k}$ with $k\ge2$ and
$\7m\subset\FF[T_{1},\dots,T_{n}]$ a maximal ideal.

\KAP{Affinely}{Affinely homogeneous surfaces}  

For the rest of the paper we consider only nilpotent algebras that are
{\sl commutative.} The algebras of finite dimension of this type are
precisely the maximal ideals of commutative Artinian local
algebras. 

It is well-known that the graph $S:=\{(x,t)\in V\oplus\FF:t=\6q(x)\}$
of every quadratic form $\6q$ on a vector space $V$ is affinely
homogeneous.  On the other hand, for given vector space $W$ of finite
dimension, the vast majority of smooth algebraic hypersurfaces in $W$
of degree $>2$ is not affinely homogeneous. In fact, it is not easy to
find an affinely homogeneous hypersurface of higher degree at all.  In
this section we show that the nil-surfaces associated with a certain
class of commutative nilpotent algebras $\5N$ are affinely homogeneous
varieties and can have arbitrary high degrees. More precisely, we have
a positive result in case where $\5N$ admits some sort of a
$\ZZ^{+}$-grading.

\Definition{NC} Let $\5N$ be a nilpotent algebra, $\pi$ a projection
on $\5N$ and $\5N=\bigoplus_{k\in\ZZ^{+}}\5N_{k}$ a vector space
decomposition. This decomposition is called \par\noindent $\bullet$ a
{\sl grading} if
$$\5N_{j}\5N_{k}\subset\5N_{j+k}\Steil{for all}j,k>0\;,$$ 
$\bullet$
 a {\sl $\pi$-grading} if
$$\5N_{j}\5N\;\subset\;\;\bigoplus_{\ell>j}\5N_{\ell}\Steil{for
all}j>0\quad\hbox{and}$$
$$\pi(\5N_{j_{1}}\5N_{j_{2}}\!\cdots\5N_{j_{r}})\;\subset\;
\pi(\5N_{j_{1}+j_{2}+\dots+j_{r}})$$ holds for every finite sequence
$j_{1},j_{2},\dots,j_{r}$ in $\ZZ^{+}$.

\medskip A quite special sort of grading is what usually is called a
{\sl canonical grading:} To every nilpotent algebra $\5N$ associate
the graded algebra
$$\gr(\5N):=\bigoplus_{k>0}\qu{\5N^{k}}{\5N^{k+1}}\,,$$ where for
every $j,\,k>0$ and every $x\in\5N^{j}$, $\,y\in\5N^{k}$ the product
of the residue classes $x{+}\5N^{j+1}$ and $\,y{+}\5N^{k+1}$ is
$\,xy{+}\5N^{j+k+1}\,$. It is quite rare that $\5N$ and $\gr(\5N)$ are
isomorphic as algebras. But if there exists an algebra isomorphism
$\phi:\5N\to\gr(\5N)$, the grading $\5N=\bigoplus\5N_{k}$ with
$\5N_{k}=\phi^{-1}(\qu{\5N^{k}}{\5N^{k+1}})$ is called a canonical
grading. Clearly $\5N_{k}=0$ for all $k>\nu(\5N)$ in this case. 

Gradings and $\id$-gradings on $\5N$ coincide.  Graded nilpotent
commutative algebras exist for every nil-index $\nu\,$ -- for instance
the maximal ideal of $\FF[X]/(X^{\nu+1})$ has an obvious canonical
grading. On the other hand, not every nilpotent commutative algebra
has a grading, see Section \ruf{counter} for counterexamples. In
general, a gradable nilpotent algebra $\5N$ may not have a grading
with $\5N_{1}\ne0$. A simple example with this phenomenon is the
commutative algebra $\5N=\FF x\oplus\FF y\oplus\FF x^{2}\oplus\FF
x^{3}$ with generators $x,y$ satisfying
$x^{4}=y^{3}=xy=x^{3}-y^{2}=0$. Then we get a grading of $\5N$ if we
denote the summands successively by $\5N_{2}$, $\5N_{3}$, $\5N_{4}$,
$\5N_{6}$. Notice that this $\5N$ does not admit a canonical
grading. Indeed, the annihilators of $\5N$ and $\gr(\5N)$ have
dimension 1 and 2 respectively.

\bigskip What about affine homogeneity of $S_\pi$ for arbitrary
commutative nilpotent algebras $\5N$ and $\5N_{[1]}$-ranged projections
$\pi$?  For some time the answer of this question was beyond our reach
as all our attempts to prove it for general commutative nilpotent
algebras (even when restricted to the special case $\dim(\5N_{[1]})=1$)
failed in case of nil-index $\ge 5$.  However, contrary to the
expectation (expressed as conjecture in \Lit{ISAA}, [\ruf{ISAE}a])
counterexamples do exist.  Anticipating the answer, which will be
extensively discussed in Section \ruf{counter}, we have:

{\parindent=10pt\item{$\bullet$} There exist commutative nilpotent
algebras $\5N$, such that $S_\pi$ is not affinely homogeneous (any
such algebra cannot be graded).  \item{$\bullet$} There exist
commutative nilpotent algebras $\5N$ without a grading but still with
affinely homogeneous $S_\pi$.\par}

\medskip Now we resume our investigation by proving the main result of
this section. For every pair of vector spaces $V,W$ the affine group
$\Aff(V)$ acts from the right on the space of all polynomial mappings
$f:V\to W$ and we denote by \smallskip\line{\hfill
$\6A_{f}:=\{g\in\Aff(V):f\circ g^{-1}=f\}$\hfill} \noindent the
isotropy subgroup at $f$. Clearly, $\6A_{f}$ leaves every level set
$f^{-1}(c)$, $\,c\in f(V)$, invariant.

\Theorem{NB} Let $\5N$ be a commutative nilpotent algebra and $\pi$ a
projection on $\5N$. Assume that $\5N=\bigoplus\5N_{k}$ is a
$\pi$-grading and that there exists an integer $d>0$ with
$\pi(\5N)\subset\5N_{d}$ and $\pi(\5N_{k})=0$ for all $k\ne d$. Then
for $\6f:=\pi\circ\exp_{1}$ the affine subgroup
$\6A_{\6f}\subset\Aff(\5N)$ acts transitively on
$S=\6f^{-1}(0)$.\nline In case $\pi(\5N)\subset\Ann(\5N)$ the group
$\6A_{\6f}$ even acts transitively on every level set
$\6f^{-1}(c)=S+c$ with $c\in\pi(\5N)$.

\Proof Fix an arbitrary point $a\in S$. For every $j\ge 0$ consider
the following condition:
$$\6A_{\6f}(a)\cap\bigoplus_{k>j}\5N_{k}\ne\emptyset\;.\leqno{(\star)}$$
\vskip-10pt\noindent It is clear that for the first claim in the
Theorem we only have to show that $(\star)$ holds for all $j$ since
then $0\in\6A_{\6f}(a)$, or equivalently $a\in\6A_{\6f}(0)$.  \nline
We first show by induction over $j$ that $(\star)$ is valid for every
$j<d\,$: For $j=0$ nothing has to be shown. Now fix an arbitrary
integer $j$ with $0<j<d$. As induction hypothesis we then may  assume
$a\in\bigoplus_{k\ge j}\5N_{k}$.  Set $\5N_{0}:=\FF\cd\One$. Every
$x\in\5N^{\,0}$ has a unique decomposition $x=x_{0}+x_{1}+\dots$ with
$x_{k}\in\5N_{k}$ for all $k\ge0$. In particular,
$a=a_{j}+a_{j+1}+\dots\;$ with $a_{k}\in\5N_{k}$ for all $k\ge j$. We
trivially extend $\6f$ to a function $\hat{\6f}$ on $\5N^{\,0}$, more
precisely, $\hat{\6f}(s\One+x):=\6f(x)$ for all $s\in\FF$ and
$x\in\5N$.\nline Denote by $\5F$ the space of all polynomial maps
$\5N^{\,0}\to\5N_{d}$ of degree $\le d$. Then $\hat{\6f}\in\5F$. We
identify via $x\leftrightarrow\One{+}x$ the nilpotent algebra $\5N$
with the affine hyperplane $\5U:=\One{+}\5N$ in $\5N^{\,0}$ and
$\6A_{\6f}$ with the subgroup
$$\6G:=\{g\in\GL(\5N^{\,0}):g(\5U)=\5U\steil{and}\hat{\6f}(gx)=
\hat{\6f}(x)\steil{for all}x\in\5U\}\;.$$ With $a_{j}\in\5N_{j}$ from
above define $\lambda=\lambda_{j}\in\End(\5N^{\,0})$ by
$$\lambda(x):=\;\;a_{j}\Big(-x_{0}+{1\over d-j}\sum_{k=1}^{d-j} k\,
x_{k}\Big)\;.\leqno{(\star\star)}$$ Then $\lambda(\One)=-a_{j}\in\5N$
and $\lambda(\5N)\subset\bigoplus_{k>j}\5N_{k}$. For
$g:=\exp(\lambda)\in\GL(\5N^{\,0})$ ($\lambda$ is nilpotent) we
therefore have $g(\One{+}a)=\One{+}b$ for some $b$ in
$\bigoplus_{k>j}\5N_{k}$. It is enough to show $g\in\6G$ since then
$b\in\6A_{\6f}(a)$ by the above identifications. The identity
$g(\5U)=\5U$ is obvious. It remains to compute $\hat{\6f}\circ g$ on
$\5U$.  This can be done in terms of the following nilpotent operator
$\xi\in\End(\5F)\,$:
$$\xi(f)x:=f'(x)(\lambda x)\Steil{for
all}f\in\5F\steil{and}x\in\5N^{\,0}\;,$$ where
$f'(x)\in\Hom(\5N^{\,0},\5N_{d})$ is the formal derivative of $f$ at
$x$. From $\lambda(\5N^{\,0})\subset\5N$ we conclude that $\xi(f)$
vanishes on $\5U$ as soon as $f$ has the same property.  For all
$f\in\5F$ we have the generalized Taylor's formula
$$f\circ\exp(\lambda)=\exp(\xi)(f)=f+
\xi(f)+{1\over2}\xi^{2}(f)+\dots\;\;.$$ It therefore remains to show
that $\xi(\hat{\6f})$ vanishes on $\5U$. Now for every $x\in\5U$ we
have $x_{0}=\One$ and
$$\xi(\hat{\6f})x=\pi\sum c_{\nu}x_{1}^{(\nu_{1})}x_{2}^{(\nu_{2})}\cdots
x_{d}^{(\nu_{d})}\,,$$ where the sum is taken over all multi indices
$\nu\in\NN^{d}$ with $\nu_{1}+2\nu_{2}+\dots+d\nu_{d}=d-j$, and
$c_{\nu}\in\5N_{j}$ are certain factors. Fix such a multi index
$\nu$. For simpler notation we put $x^{(-1)}:=0$ for every
$x\in\5N^{\,0}$.  Then we have {\klein
$$\eqalign{c_{\nu}^{}\,x_{1}^{(\nu_{1})}x_{2}^{(\nu_{2})}\cdots
x_{d}^{(\nu_{d})}&=\sum_{k=1}^{d-j}{k\over d-j}a_{j}
x_{k}\dd{x_{j+k}}\Big(x_{1}^{(\nu_{1})}\cdots
x_{k}^{(\nu_{k}-1)}\cdots x_{j+k}^{(\nu_{j+k}+1)}\cdots
x_{d}^{(\nu_{d})}\Big) \cr
&\phantom{XXXX}-a_{j}\!\dd{x_{j}}\Big(x_{1}^{(\nu_{1})}\cdots
x_{j}^{(\nu_{j}+1)}\cdots x_{d}^{(\nu_{d})}\Big) \cr
&=\Big(\sum_{k=1}^{d-j}{k\nu_{k}\over
d-j}-1\Big)a_{j}x_{1}^{(\nu_{1})} x_{2}^{(\nu_{2})}\cdots
x_{d}^{(\nu_{d})}\;=\;0 }$$ }\noindent since $\nu_{k}^{}=0$ for
$k>d-j$. This proves the induction step and hence $(\star)$ for all
$j<d$, that is, we may assume $a\in\bigoplus_{k\ge d}\5N_{k}$. \nline To
finish the proof of the first statement we notice that
$\pi\big((x+b)^{k}\big)=\pi(x^{k})$ holds for all $x\in\5N$, $k\ge1$
and every $b\in\big(\id-\pi\big)\big(\bigoplus_{k\ge d}\5N_k\big)$. This
implies that the translation $x\mapsto x+b$ for every such $b$ belongs
to $\6A_{\6f}$. As a consequence we may assume $a\in\pi(\5N_{d})$.
But $a$ is also in $S$ by assumption and $S\cap\pi(\5N_{d})=\{0\}$,
that is, $a=0$.

\smallskip\noindent Now suppose $\pi(\5N)\subset\Ann(\5N)$ and let
$\6B$ be the subgroup of all $g\in\6A_{\6f}$ that commute with all
translations $x\mapsto x+c\,$, $\,c\in\pi(\5N_{d})$. In every
induction step above the operator $\lambda$ vanishes on
$\pi(\5N)\subset\5N^{\,0}$. This implies that $\6B$ is already
transitive on $S$ and the second claim follows. \qed

\smallskip In the proof of \ruf{NB} we have identified $\5N$ with
$\;\One+\5N$ via the identification $x\leftrightarrow\One{+}x$. In
case $\5N_{d}$ is the annihilator $\5N_{[1]}$ of $\5N$ in \ruf{NB}, the
operator $\lambda$ in $(\star\star)$ corresponds via the
identification to the affine transformation $T:\5N\to\5N$, where
$$T={1\over(d-j)}D-a_{j}\steil{with}D\in\der(\5N)\steil{defined
by}D(x)=a_{j}\sum_{k>0}k\, x_{k}\;.\Leqno{ABC}$$ The proof in \Lit{ISAE} for the
special case, where $\5N_{d}$ has dimension 1 and is the annihilator
of $\5N$, is also based on these nilpotent derivations $D$.

\medskip For base field $\FF\in\{\RR,\CC\}$ Theorem \ruf{NB}
essentially is already contained in \Lit{FEKA}, see also \Lit{FIKK}
for a special version with $\FF=\CC$. For the special case that
$\pi(\5N)$ is the annihilator of $\5N$ and this annihilator has
dimension 1 see also \Lit{ISAE}.

\bigskip In Theorem \ruf{NB} the group $\6A_{\6f}$ is not the full
affine group $\Aff(S)$. Indeed,
$$\theta_{t}:=\bigoplus_{k>0}t^{k}\id_{|\5N_{k}}\;\in\;\GL(S)\Leqno{ZQ}$$
satisfies $\6f\circ\theta_{t}=t^{d}\6f$ for every $t\in\FF^{*}$. As a
consequence, if $\pi$ has 1-dimensional range in
$\5N_{d}\cap\Ann(\5N)$, the group $\Aff(S)$ has at most $d$ orbits in
$\5N$. In particular, in case $\;\FF=\CC\;$ this group has only two
orbits in $\5N$, the hypersurface $S$ and its open connected
complement $\5N\backslash S$.  In case $\;\FF=\RR\;$ the connected
identity component $\Aff(S)^{0}$ has three orbits in $\5N$, the
hypersurface $S$ and both sides of the complement $\5N\backslash S$.

In case $\5N=\bigoplus\5N_{k}$ is a grading in \ruf{NB}, the operator
$\lambda:=\bigoplus_{k>0}k\id_{|\5N_{k}}\in\der(\5N)$ is
diagonalizable over $\FF$ and has only positive integers as
eigen-values. Conversely, if $\5N$ is an arbitrary (commutative)
nilpotent algebra and $\lambda\in\der(\5N)$ is diagonalizable over
$\FF$ with spectrum in $\ZZ^{+}$, then $\5N=\bigoplus\5N_{k}$ is a
grading, where $\5N_{k}$ for every $k$ is the $k$-eigenspace of
$\lambda$.

\bigskip As already mentioned, not every nilpotent algebra $\5N$ has a
grading, compare Section \ruf{counter} for counterexamples. But there
exists always a decomposition
$$\5N=\bigoplus_{k\in\ZZ^{+}}\5N_{k}\,,\Steil{with}
\5N_{j}\5N_{k}\subset\bigoplus_{\ell\ge j+k}\5N_{\ell}\steil{for
all}j,k>0\,.\Leqno{CN}$$ Indeed, choose $\5N_{k}$ in such a way that
$\5N^{k}=\5N_{k}\oplus\5N^{k+1}$ for all $k>0$. In general, a
non-gradable nilpotent algebra may have a $\pi$-grading with
non-trivial projection $\pi$.  As a trivial example for this
phenomenon choose for fixed non-gradable $\5N$ a decomposition
\Ruf{CN} and let $\pi$ on $\5N$ be the canonical projection onto
$\5N_{2}$.

\medskip The following result gives a lower bound for the size of the
0-orbit under the affine group $\Aff(S)$ in case $\pi$ has range in
the annihilator of $\5N$: Fix $k>0$ and consider the ideal $\5N_{[k]}$
as defined in $\ruf{NA}$. Then for $\6f:=\pi\circ\exp_{1}$ and
$S:=\6f^{-1}(0)$ as before the intersection $S\cap\5N_{[k]}$ is a
smooth subvariety with dimension $\dim(\5N_{[k]}/\5N_{[1]})$. For
every $a\in S\cap\5N_{[3]}$ and $\rho:=(\id+\pi)/2\in\End(\5N)$ define
the affine transformation $h_{a}$ on $\5N$ by
 $$h_{a}(x):=x-\rho(ax)+a\;.$$
 
\Proposition{ZT} Let $\5N$ be an arbitrary commutative nilpotent
algebra, let $\pi$ be a projection on $\5N$ with range in the
annihilator of $\5N$ and let $S:=S_{\pi}$.  Then $\{h_{a}:a\in
S\cap\5N_{[3]}\}$ is contained in $\6A_{\6f}$ and generates a subgroup
acting transitively on $S\cap\5N_{[3]}$. In particular, $S$ is
affinely homogeneous if $\5N$ has nil-index $\le3$.

\Proof Fix $a\in S\cap\5N_{[3]}$. Then $h_{a}\in\Aff(\5N)$ since the
operator $x\to\rho(ax)$ is nilpotent on $\5N$.  A simple computation
shows $\6f\circ h_{a}=\6f$ and also that $h^{}_{a},h_{a}^{-1}$ leave
$\5N_{[3]}$ invariant. The first claim follows with $h_{a}(0)=a$. The
second follows from $\5N_{[3]}=\5N$ in case of
$\nu(\5N)\le3$.\qed

\KAP{Admissible}{Admissible algebras}

For the rest of this paper we deal only with commutative nilpotent
algebras $\5N$ of finite dimension over $\FF$ such that the
annihilator $\5N_{[1]}$ is of dimension $1$.  For simplicity we call
such algebras {\sl admissible algebras}. These are just the maximal
ideals of Gorenstein algebras of finite vector space dimension $\ge2$
over $\FF$.

In this and the subsequent sections we construct several objects,
universally associated with a given admissible algebra $\5N$. These
will encode enough information to characterize the admissible algebra
up to isomorphy.  Roughly speaking, we define a certain family
$\Sigma$ of smooth hypersurfaces $S\subset \5N$ such that each of its
members determines $\5N$. We also establish a natural duality between
the points of a given hypersurface $S\in\Sigma$ and the members of
$\Sigma$ itself. In the next section we construct a set of
$\FF$-valued polynomials, so-called nil-polynomials, closely related
to the hypersurfaces $S\in\Sigma$. We also determine how the algebra
structure of $\5N$ can be reconstructed from an associated
nil-polynomial $\6p$ (in fact the knowledge of the quadratic and cubic
terms of $\6p$ turns out to be sufficient.)

We start with some preparations.  We call every projection
$\pi=\pi^{2}\in\End(\5N)$ with range $\pi(\5N)=\5N_{[1]}$ an {\sl
admissible projection} on $\5N$ and denote by
$\Pi(\5N)\subset\End(\5N)$ the subvariety of all admissible
projections. Every $\pi\in\Pi(\5N)$ is uniquely determined by its
kernel $\5K=\ker(\pi)$ that satisfies $\5N=\5K\oplus\5N_{[1]}$.
Further, every projection $\pi\in \Pi(\5N)$ gives rise to the
algebraic smooth hypersurface, compare \Ruf{ND},
$$ S_\pi=\{x\in\5N:\pi\circ \exp_{1}(x)=0\}\;=\;\log_{1}(\ker\pi)$$
that we also call a {\sl nil-hypersurface}. We denote by
$\Sigma(\5N):=\{S_\pi:\pi \in \Pi(\5N)\}$ the set of all such
hypersurfaces. Note that $\{0\}$ is the intersection of all $S\in
\Sigma(\5N)$.  The canonical map $\beta:\Pi(\5N)\to \Sigma(\5N)$,
$\pi\mapsto S_\pi$, is clearly surjective. Later on (see \ruf{PS}) we
will show that $\beta$ is even bijective.

All the key objects associated with $\5N$, such as the bilinear forms
$\6b_\pi:\5N\times \5N\to \5N_{[1]}$, the polynomial maps
$\6f_\pi:\5N\to \5N_{[1]}$ and the subvarieties $S_\pi\subset \5N$
depend on the choice of the projection $\pi$.  In this section we show
that in the admissible case the `essential' properties of $\6b_\pi,
\6f_\pi$ and $S_\pi$ do not depend on the projection and can be
considered as invariants associated to $\5N$ only. Further we prove
that every smooth hypersurface $S_{\pi}$ determines the admissible
algebra $\5N$ up to isomorphy.

For every $e\in\5N^{\,0}$ define the multiplication operator
$M_{e}\in\End(\5N)$ by $M_{e}(x)=ex$. Recall that $\exp:\5N\to
\5U=\One+\5N$ is a group isomorphism with inverse $\log$.

\Lemma{PA} $\Pi(\5N)$ is an affine plane of dimension
$\dim(\5N/\5N_{[1]})$ in $\End(\5N)$. In fact
$$\Pi(\5N)=\{\pi\in\
\End(\5N):\pi=\rho\circ\pi,\,\rho=\pi\circ\rho\}\Steil{for
every}\rho\in\Pi(\5N)$$ and $ \5N/\5N_{[1]} \times
\Pi(\5N)\,\longrightarrow\, \Pi(\5N)\,,\;
(x+\5N_{[1]}\,,\,\pi)\longmapsto \pi\circ M_{\exp x}\,,$ yields a
well-defined simply transitive action of the vector group
$\5N/\5N_{[1]}$ on $\Pi(\5N)\,$.

\Proof Clearly, the mapping is well defined, has values in $\Pi(\5N)$
and gives an action of $\5N/\5N_{[1]}$. The action is also free --
indeed, suppose that $\pi\circ M_{\exp b}=\pi$ for some
$\pi\in\Pi(\5N)$, $b\in\5N$. For $c:=\exp_{1}(b)$ then $\pi\circ
M_{c}=0$, that is, $\6b_\pi(c,\5N)=0$ and thus $c\in\5N_{[1]}$, see
\Ruf{NF}. But then also $b=c\in\5N_{[1]}$. \nline The action is also
transitive -- indeed, fix arbitrary $\pi,\rho\in\Pi(\5N)$. Then
$\lambda:=\rho-\pi$ vanishes on $\5N_{[1]}$ and satisfies
$\lambda=\pi\circ\lambda$.  Hence, again by the non-degeneracy of
$\6b_\pi$ on $\ker(\pi)$, we conclude $\lambda=\pi\circ M_{b}$ for
some $b\in\ker\pi$. This implies $\rho=\pi+\pi\circ M_{b}=\pi\circ
M_{\one+b}=\pi\circ M_{\exp c}$ for $c:=\log(\One+b)\in\5N$.\qed

The algebra automorphism group $\Aut(\5N)$ acts on $\Pi(\5N)$ by
conjugation, that is, by $L(\pi):=L\circ\pi\circ L^{-1}$ for all
$L\in\Aut(\5N)$ and $\pi\in\Pi(\5N)$. Then $L(\ker\pi)=\ker L(\pi)$ is
obvious.  The group $\Aut(\5N)$ also acts on $\Sigma(\5N)$ in the
obvious way and satisfies $L(S_{\pi})=S_{L(\pi)}$ for all
$L\in\Aut(\5N)$ and $\pi\in\Pi(\5N)$.

The following result generalizes Propositions 2.2. and 2.3 in
\Lit{FIKK}.

\Theorem{AB} Let $\5N,\tilde {\5N}$ be arbitrary admissible algebras
 having not necessarily the same dimension. Also let $\pi\in
 \Pi(\5N)$, $\tilde \pi\in \Pi(\tilde {\5N})$ be arbitrary admissible
 projections. Then for $S:=S_\pi$, $\,\tilde S:=S_{\tpi}$ and for
 every linear map $L:\5N\to\tilde\5N$ the following conditions are
 equivalent, provided $\dim(\tilde{\5N})>1$. \0 $L$ is an algebra
 isomorphism.  \1 $\tilde S=L(S-c)$ for some $c\in\5N$.

\noindent Furthermore, the point $c=c_{L,\pi,\tpi}$ in {\rm (ii)} is
uniquely determined by $L,\pi,\tilde\pi$ and coincides with the unique
element in $S$ satisfying $\tilde\pi=L\circ(\pi\circ M_{\exp c})\circ
L^{-1}$. Finally $$
S=\big\{c_{L,\pi,\rho}:\rho\in\Pi\big(\tilde\5N\big)\,\big\}\leqno{(*)}$$
holds for every algebra isomorphism $L:\5N\to\tilde\5N$.

\Proof \To12 ~Assume (i).  By Lemma $\ruf{PA}$ there exists $c\in\5N$
with $L^{-1}\circ\tilde\pi\circ L=\pi\circ M_{\exp c}\,$. Since
$\exp_1(c+a)=\exp_1(c)+a$ for every $a\in \5N_{[1]}$, we can assume
$\pi(\exp_{1} c)=0$, that is, $c\in S$. Then $\pi(\exp_{1} c)=0$
implies
$$\tilde\pi(\exp_{1}L(x))=\tilde\pi\circ L(\exp_{1}x)=L\circ \pi((\exp
c)(\exp_{1}x))=L\circ\pi(\exp_{1}(x+c))$$ for all $x\in\5N$, that is,
$L(x)\in\tilde S$ if and only if $x+c\in S$.

\noindent\To21 ~Since $\dim \tilde{\5N}>1$ also $\nu(\tilde\5N)>1$ and
the linear span of $\widetilde S$ is $\tilde{\5N}$.  Then $L$ is an
epimorphism and also $\nu(\5N)>1$. We show that $L$ is also injective:
Note first that since the quadratic part of $\pi\circ\exp_2$ is
non-degenerate on $\ker(\pi)$, $S$ is not invariant under any
non-trivial translation. If $\ker(L)\ne 0$ then $L^{-1}(\tilde
S)=\ker(L){+}S{-}c\,$ would be Zariski dense in $\5N$ and a proper
algebraic subset of $\5N$ at the same time, a contradiction.\nline
Since $L-L(c)$ provides an affine equivalence between $S$ and $\tilde
S$ we can use the analytic proof of Prop. 2.3 in \Lit{FIKK} to obtain
that $L$ is an $\FF$-algebra isomorphism in the special case
$\FF=\CC$. We reduce the case of a general field to this special
result by a Lefschetz principle type argument. To begin with we denote
by $\3K$ the set of all subfields $\KK\subset\FF$ that are obtained by
adjoining a finite subset of $\FF$ to the prime field of $\FF$. It is
well known that every $\KK\in\3K$ is isomorphic to a subfield of
$\CC$.\nline Now let $e\in\5N$ be an arbitrary but fixed element. Then
it is enough to show that $L(e^{2})=L(e)^{2}$: Choose a linear basis $
B$ of $\5N$ containing a basis of $\pi^{-1}(0)$ and a basis of
$\5N_{[1]}$.  Then there exists a field $\KK\in\3K$ such that the
linear span $\5B:=\langle B\rangle_{\KK}$ contains $e$ and is a
$\KK$-subalgebra of $\5N$.  By the choice of $ B$ the intersection
$\5B\cap\5N_{[1]}$ has dimension 1 over $\KK$ and is the annihilator
of $\5B$. Also, $S\cap\5B$ is a smooth hypersurface over $\KK$ in
$\5B$. In the same way choose a linear basis $\tilde B$ of $\tilde\5N$
containing a basis of $\tilde\pi^{-1}(0)$ and a basis of
$\tilde\5N_{[1]}$. Adjoining a suitable finite subset of $\FF$ to
$\KK$ we may assume in addition without loss of generality that
$\tilde\5B:=\langle\tilde B\rangle_{\KK}$ contains $c$ and also is a
$\KK$-subalgebra of $\tilde\5N$. Enlarging $\KK$ again within $\3K$ if
necessary, we may even assume that the affine transformation
$A:=L-L(c)$ maps $\5B$ onto $\tilde\5B$. Clearly $A$ maps $S\cap\5B$
onto $\tilde S\cap\tilde\5B$.  \nline We now consider $\KK$ as
subfield of $\CC$. We then get the complex nilpotent algebras
$\5B\otimes_{\KK}\CC$ and $\tilde\5B\otimes_{\KK}\CC$ with
annihilators $(\5B\cap\5N_{[1]})\otimes_{\KK}\CC$ and
$(\tilde\5B\cap\tilde\5N_{[1]})\otimes_{\KK}\CC$ respectively, each
having complex dimension 1 over $\CC$. The $\KK$-affine map $A_{|\5B}$
extends to a $\CC$-affine map $\5B\otimes_{\KK}\CC\to$
$\tilde\5B\otimes_{\KK}\CC$ sending the corresponding complex
hypersurfaces onto each other. By Proposition 2.3 in \Lit{FIKK} then
$M_{|\5B}\otimes_{\KK}\id_{\CC}$ is an algebra isomorphism, implying
$L(e^{2})=L(e)^{2}$. Together with the first step this proves \TO12.

\noindent Next, assume that $c$ in {\rm (ii)} is not uniquely
determined. Then there exists $a\in\5N$ with $a\ne0$ and $S=S+a$. For
$\5K:=\pi^{-1}(0)$ we have $\5N=\5K\oplus\5N_{[1]}$ and
$S=\{(y,f(y)):y\in\5K\}$ is the graph of the polynomial map
$f:\5K\to\5N_{[1]}$ given by $f(y)=-\pi(\exp_{2}(y))$. In particular,
$f=(b,f(b))$ for some $b\in\5K$ and
$$f(y+tb)=f(y)+f(tb)\Steil{for all}y\in\5K\leqno{(**)}$$ and all
$t\in\ZZ$.  Since $f$ is a polynomial map $(**)$ even holds for all
$t\in\FF$.  Comparing terms that are linear in $y$ as well as in $t$
we get $\pi(b y)=0$ for all $y\in\5K$. But the quadratic form
$\pi(y^{2})$ is non-degenerate on $\5K$, implying $b=0$ in
contradiction to $f\ne0$.

\noindent Finally, for the proof of $(*)$ we may assume
$\tilde\5N=\5N$ and $L=\id_{\5N}$. Fix an arbitrary $c\in S$ and put
$e:=\exp_{1}(c)$. Then $\pi(e)=0$ and $\rho:=\pi+\pi\circ M_{e}$ is an
admissible projection on $\5N$. This implies $S_{\rho}=S-c$ and
consequently $c=c_{L,\pi,\rho}$.\qed
 
\Corollary{AC} The algebras $\5N$, $\tilde\5N$ are isomorphic if and
only if $S$, $\tilde S$ are affinely equivalent.\Formend

\Corollary{AE} Under the same assumptions as in \ruf{AB}, for every
linear map $L:\5N\to\tild\5N$ the following conditions are
equivalent. \smallskip \0 $\tilde S=L(S)$.\1 $L$ is an algebra
isomorphism with $\tilde\pi=L(\pi)\;\;\;(=L\circ\pi\circ L^{-1})$.

\Proof Assume (i). Then $L$ is an algebra isomorphism with
$\tilde\pi=L(\pi\circ M_{\exp c})=L(\pi)$ for $c=0$ by Theorem
\ruf{AB}. ~The converse implication is trivial.\qed

\Corollary{} $\Aut(\5N)\cap\Aff(S)=\GL(S)$.\Formend

\medskip Next we show equivalences between the various sets. In
particular, for every fixed $S\in \Sigma(\5N)$, we give a duality
between points in $S$ and surfaces in $\Sigma$ itself.

\Lemma{PS} {\bf (Duality)} Let $\5N$ be an admissible algebra and
$\pi\in\Pi(\5N)$. Then the mappings
$$\diagram{\alpha_\pi:S_\pi&\longrightarrow& \Pi(\5N)\cr
\qquad s&\longmapsto &\pi\circ M_{\exp s}},\qquad
\diagram{\beta:\Pi(\5N)&\longrightarrow& \Sigma(\5N)\cr
\qquad \rho&\longmapsto &S_\rho}\qquad\qquad
$$ are bijective and satisfy \quad $\beta\circ
\alpha_\pi(s)=S_\pi-s$\quad for all $s\in S_{\pi}$.

\Proof Let $\rho=\pi\circ M_{\exp s}$ with $s\in S_{\pi}$. Then
$x\in S_{\rho}$ is equivalent to
$$0=\pi\big(\exp(s)\exp_{1}(x)\big)=\pi\big(\exp_{1}(x+s)-
\exp_{1}(s)\big)$$ and hence to $(x+s)\in S_{\pi}$ since
$\pi(\exp_{1}(s))=0$. Bijectivity of $\alpha_\pi$ follows from the
proof of \ruf{PA} and the fact that $S_{\pi}$ as graph has the
following property: Every $x\in\5N$ is a unique sum $x=a+b$ with
$a\in\5N_{[1]}$ and $b\in\S_{\pi}$. Surjectivity of $\beta$ holds by
definition and injectivity follows from \ruf{AE}.\qed

Our next goal is to investigate the behavior of the above sets under
algebra isomorphisms. Obviously, algebra isomorphisms are functorial
in the sense that for every such isomorphism $L:\5N\to \tilde {\5N}$
one has the following well-defined maps (also denoted by the same
letter $L$):
$$\diagram{L:\Sigma(\5N)&\longrightarrow& \Sigma(\tilde{\5N})\cr
\qquad S&\longmapsto &L(S)},\qquad
\diagram{L:\Pi(\5N)&\longrightarrow& \Pi(\tilde{\5N})\cr \qquad
\pi&\longmapsto &L\circ \pi \circ L^{-1}}\qquad \Leqno{PX}$$\vskip-20pt
$$\rlap{\hskip-25pt with}L(S_\pi)=S_{L(\pi)}\;.\Leqno{LS}$$ In
particular, we have in case of $\5N=\tilde{\5N}$ the group action of
the algebra automorphism group $\Aut(\5N)\subset\GL(\5N)$ on the
affine plane $\Pi(\5N)$ by conjugation, that is, by
$L(\pi)=L\circ\pi\circ L^{-1}$ for all $L\in\Aut(\5N)$ and
$\pi\in\Pi(\5N)$. Also the affine group $\Aff(S)$ acts canonically on
the hypersurface $S\in\Sigma(\5N)$ and we show next that both group
actions are equivariantly equivalent, more precisely, define the
following map
$$\gamma:\Aut(\5N)\longrightarrow \Aff(\5N),\qquad L\mapsto
\gamma(L):=L-L(c_{L,\pi,\pi})\;,$$ with $c_{L,\pi,\pi}\in S$ as in
Theorem \ruf{AB}, see also the first part of its proof.

\Proposition{NI} Let $\5N$ be an admissible algebra, $\pi\in\Pi(\5N)$,
$S:=S_\pi$ and $\gamma$ as above. Then $\gamma$ induces a group
isomorphism $\Aut(\5N)\to \Aff(S)$. Furthermore, the diagram

$$ \diagram{\Aut(\5N)&\;\times\; &\Pi(\5N)&\longrightarrow &
\Pi(\5N)\cr \Big\downarrow\rlap{$\gamma$}
&&\Big\downarrow\rlap{$\alpha_\pi^{-1}$}
&&\Big\downarrow\rlap{$\alpha_\pi^{-1}$}\cr \Aff(S)&\;\times\;
&S&\longrightarrow & S }$$ commutes and has bijective vertical arrows,
while the horizontal arrows represent the respective group actions.

\Proof For $A:=\gamma(L)$ we have $A(S)=L(S-c_{L,\pi,\pi})=S\,$ by
(ii) of \ruf{AB}, i.e., $\gamma$ yields a map $\Aut(\5N)\to \Aff(S)$,
which by Theorem \ruf{AB} is a bijection onto $\Aff(S)\,$.  The
inverse of $\gamma$ is just the mapping that associates to every
$A\in\Aff(S)$ its linear part $L:=A-A(0).$ This implies that $\gamma$
is a group isomorphism.  The commutativity of the diagram can be seen
as follows: Direct consequence of \ruf{LS} is the commutativity of the
diagram
$$ \diagram{\Aut(\5N)&\;\times\; &\Sigma(\5N)&\longrightarrow &
\Sigma(\5N)\cr \big\Vert &&\big\uparrow\rlap{$\beta$}
&&\big\uparrow\rlap{$\beta$}\cr \Aut(\5N)&\;\times\;
&\Pi(\5N)&\longrightarrow & \Pi(\5N)}$$ with bijective vertical
maps. Hence, it suffices to prove the commutativity of
$$ \diagram{\Aut(\5N)&\;\times\; &\Sigma(\5N)&\longrightarrow &
\Sigma(\5N)\cr \big\downarrow\rlap{$\gamma$}
&&\big\uparrow\rlap{$\beta\circ \alpha_\pi$}
&&\big\uparrow\rlap{$\beta\circ\alpha_\pi$}\cr \Aff(S)&\;\times\;
&S&\longrightarrow & S\rlap{\qquad.} }$$
According to \ruf{PS} and \ruf{AB}
we have for arbitrary $L\in \Aut(\5N)$ and $t\in S=S_\pi$
$$ \diagram{\rlap{$\hskip-29pt\big($}L&\;,\;
&S-t\,\big)&~~\longrightarrow & ~~~L(S)-L(t)\rlap{
$=S+L(c_{L,\pi,\pi})-L(t)$}\cr \big\downarrow\rlap{$\gamma$}
&&\big\uparrow\rlap{$\beta\circ \alpha_\pi$}
&&\big\uparrow\rlap{$\beta\circ\alpha_\pi$}\cr\big(L-L(c_{L,\pi,\pi})&\;,\;
&t\rlap{\hskip9pt$\big)$}&~~\longrightarrow &
\rlap{\hskip-20pt$L(t)-L(c_{L,\pi,\pi})$\qquad.} } \mskip140mu$$
\vskip-24pt\qed

\bigskip \noindent{\bf Remark.} Proposition \ruf{NI} remains true also
in the trivial case $\dim(\5N)=1$, although then $S$ consists of a
single point and $\Aut(\5N)=\GL(\5N)$ holds. Indeed, by our definition
in this case $\Aff(S)=\GL(\5N)$ as well. In all other cases, $S$ is
total in $\5N$ and $\Aff(S)$ acts effectively on $S$. The proposition
implies that the orbit structure for $\Aff(S)$ in $S$ is isomorphic to
the orbit structure for $\Aut(\5N)$ in $\Pi(\5N)$. In particular, both
group actions of $\Aut(\5N)$ and $\Aff(S)$ are transitive as soon as
one of it has this property. For the special case $\FF=\RR$ or
$\FF=\CC$ this is essentially the content of Theorem 2.3 in
\Lit{ISAE}: There the space $\TT$ of all hyperplanes in $\5N$
transversal to $\5N_{[1]}$ is introduced, which via
$\pi\leftrightarrow\ker(\pi)$ can be canonically identified with our
space $\Pi(\5N)$. In addition a certain subgroup
$G_{\pi}\subset\Aff(S)$ is introduced, and as Theorem 2.2 in
\Lit{ISAE} it is proved that $\Aut(\5N)$ acts transitively on $\TT$ if
and only if $G_{\pi}$ acts transitively on $S$. Then Theorem 2.3 says
$G_{\pi}=\Aff(S)$ in case of base field $\RR$ or $\CC$.

\medskip Theorem \ruf{AB} together with Proposition \ruf{NI} implies
the following result.

\Corollary{NE} For every admissible algebra $\5N$ the following
conditions are equivalent. \0 For some (and hence every)
$\pi\in\Pi(\5N)$ the hypersurface $S_{\pi}$ is affinely
homogeneous. \1 For all $\pi,\tilde\pi\in\Pi(\5N)$ the hypersurfaces
$S_{\pi},S_{\tpi}$ are linearly equivalent. \1 The group $\Aut(\5N)$
acts transitively on $\,\Pi(\5N)$. \1 The group $\Aut(\5N)$
acts transitively on $\,\Sigma(\5N)$.

\Proof \To12~ Let $\pi$, $\tilde\pi$ be admissible projections and
assume that $S_{\pi}$ is affinely homogeneous. Then
$S_{\tpi}=S_{\pi}{-}c$ for some $c\in S$ by Theorem \ruf{AB}.
For every $A\in\Aff(S_{\pi})$ with $A(0)=c$, the linear
transformation $x\mapsto\alpha(x){-}c$ maps $S_{\pi}$ onto $S_{\tpi}$.

\noindent\To21~ Assume (ii) and fix $\pi\in\Pi(\5N)$ together with an
arbitrary point $c\in S_{\pi}$. By Theorem \ruf{AB} there exists
$\tilde\pi\in\Pi(\5N)$ with $S_{\tpi}=S_{\pi}{-}c$. By assumption
there exists $g\in\GL(\5N)$ with $g(S_{\pi})=S_{\tpi}$. The
transformation $x\mapsto g(x){+}c$ is in $\Aff(S_{\pi})$ and maps the
origin to $c$.

\noindent\TO13 This follows immediately from Proposition
\ruf{NI}.\qquad\TO34 is trivial.\qed

\medskip As an immediate consequence of Theorem \ruf{NB} we state:

\Corollary{} For every graded admissible algebra $\5N$ conditions
{\rm (i) - (iv)} in \ruf{NE} hold. \qed
\Formend

\medskip Proposition \ruf{NI} says, in particular, that the groups
$\Aut(\5N)$ and $\Aff(S)$ are isomorphic. A careful inspection of the
corresponding proofs reveals that under the assumptions of Theorem
\ruf{NB} and of Proposition \ruf{ZT} these groups contain unipotent
subgroups of dimension $\dim(\5N){-}1$. In case $\5N$ has a grading,
we even get $\dim\Aut(\5N)\ge\dim(\5N)$ since then
$\theta_{s}\in\Aut(\5N)$ as defined in \Ruf{ZQ}. The same argument
gives $\dim\der(\5N)\ge\dim(\5N)=\dim(\5N^{\,0})-\dim(\5N_{[1]})$ in
the graded case, compare also Proposition 2.3 in \Lit{XUYA} in case
$\FF=\CC$. For every cyclic nilpotent algebra $\5N$ equality holds.

\medskip\noindent{\bf The infinitesimal analogon.} As shown in
\ruf{NI} the groups $\Aut(\5N)$, $\Aff(S)$ are always isomorphic.  In
case $\FF=\RR,\CC$ these groups are even isomorphic as Lie groups,
implying that then also the corresponding Lie algebras $\der(\5N)$,
$\aff(S)$ are isomorphic.  Besides $\der(\5N)$ also a Lie algebra
$\aff(S)$ can be canonically defined for arbitrary base fields, but a
priori there is no reason why these Lie algebras should be isomorphic
also in case $\FF\ne\RR,\CC$:\nline Fix an admissible algebra $\5N$
and an admissible projection $\pi$ on $\5N$. Put $S:=S_{\pi}$, that is
$S=\6f^{-1}(0)$ for $\6f:=\pi\circ\exp_{1}$.  For every $x\in S$ then
$\T_{x}(S):=\ker(\6f'(x))$ is the {\sl tangent space} at $x$, where
$\6f'(x)=\pi\circ M_{\exp x}\in\End(\5N)$ is the formal derivative of
$\6f$ at $x$ and, as defined above, $M_{y}\in\End(\5N)$ is the
multiplication operator $z\mapsto yz$.\nline Denote by $\aff(S)$ the
linear space of all affine transformations $A:\5N\to\5N$ that are
'tangent` to $S$, that is, satisfy $A(x)\in\T_{x}(S)$ for all $x\in
S$. Then $\aff(S)$ is a Lie algebra with respect to $[A,B]=A'\circ
B-B'\circ A$, where the derivative $A'=A-A(0)$ is the linear part of
$A$. A subalgebra is $\gl(S):=\gl(\5N)\cap\aff(S)$.

\Proposition{UO} For every $D\in\End(\5N)$ the following conditions are
equivalent. \0 $D\in\der(\5N)$. \1 $D-v\in\aff(S)$ for some
$v\in\5N$.

\noindent Furthermore, the vector $v=v^{}_{D,\pi}$ in {\rm (ii)} is
uniquely determined by $D$ and coincides with the unique element $v$
in $\T_{0}S=\ker \pi$ satisfying $[\pi,D]=\pi\circ M_{v}$. Also
$D\mapsto D-v^{}_{D,\pi}$ induces a Lie algebra isomorphism
$\der(\5N){\buildrel\approx\over\to}\aff(S)$.

\Proof Assume (i). Then $D(\5N_{[1]})\subset\5N_{[1]}$ for the
annihilator $\5N_{[1]}$ and hence $\pi+[\pi,D]\in\Pi(\5N)$. By Lemma
$\ruf{PA}$ there exists $c\in\5N$ with $\pi+[\pi,D]=\pi\circ M_{\exp
c}\,$, that is, $[\pi,D]=\pi\circ M_{v}$ for $v=\exp_{1}(c)$. It is no
restriction to assume $v\in\ker(\pi)=\T_{0}S$. Consider the affine
transformation $A:=D-v$ on $\5N$. Then $(\exp x)D(x)=D(\exp_{1}x)$ and
$\pi(v\exp x)=\pi(v\exp_{1} x)$ imply
$$\big(\pi\circ M_{\exp
x}\big)A(x)=\pi\circ(D-M_{v})(\exp_{1}x)=D\circ\pi(\exp_{1}x)=0$$ for
all $x\in S$. This means $A(x)\in\T_{x}S$ and (ii) is proved.\nline
Assume conversely (ii). We have to show $D(c^{2})=2cD(c)$ for all
$c\in\5N$. In case $\FF=\CC$ this follows with the affine vector field
$\xi:=(D(x)-v)\dd x$ on $\5N$ and applying Theorem \ruf{AB} to the
1-parameter subgroup $\exp(t\xi)$ of $\Aff(S)$. The case of general
base field can be reduced to $\CC$ by a Lefschetz type argument
similar to the one used in the proof of $\ruf{AB}$, we omit the
details. Also the remaining claims follow as in the proof of
\ruf{AB}.\qed

\KAP{Nil}{Nil-polynomials}

For every admissible algebra $\5N$ with annihilator $\5N_{[1]}$ we
call a linear form $\omega:\5N\to\FF$ a {\sl pointing on $\5N$} if
$\omega(\5N_{[1]})=\FF$ (in analogy to function spaces where points in
the underlying geometric space induce linear forms with certain
properties). Also, $\5N$ with a fixed pointing is called a {\sl
pointed algebra}. For every vector space $W$ of finite dimension we
denote by $\FF[W]$ the algebra of all ($\FF$-valued) polynomials on
$W$. Since in characteristic zero every field is infinite, we do not
distinguish between polynomials in $\FF[W]$ and the polynomial
functions $W\to\FF$ they induce.

\Definition{RS} $\6p\in\FF[W]$ is called a {\sl nil-polynomial}
associated to the admissible algebra $\5N$ if there exists a pointing
$\omega$ on $\5N$ and a linear isomorphism
$\phi:W\to\ker(\omega)\subset\5N$ such that
$\6p=\omega\circ\exp_{2}\circ\,\phi$.

Notice that we do not exclude the trivial case $W=0$ with
nil-polynomial $\6p=0$. Let us agree that this $\6p$ has degree
$-\infty$.

To every pointing $\omega$ on the admissible algebra $\5N$ there
exists a unique admissible projection $\pi$ on $\5N$ and a unique
linear isomorphism $\psi:\5N_{[1]}\to\FF$ with $\omega(x)=\psi(\pi x)$
for all $x\in\5N$, and conversely, every pointing on $\5N$ is obtained
this way. To $\pi$ we have associated the hypersurface
$S_{\pi}\subset\5N$, compare \Ruf{ND}.  It is easy to see that for the
nil-polynomial $\6p$ occurring in \ruf{RS} the hypersurface $S_{\pi}$
is linearly equivalent to the graph $$\Gamma_{\!\6p}:=\{(x,t)\in
W\oplus\FF:t=\6p(x)\}\,.$$

\Definition{} We say that an admissible algebra $\5N$ has {\bf
Property (AH)} if for some (and hence every) nil-polynomial $\6p$
associated to $\5N$ the graph $\Gamma_{\6p}$ is affinely homogeneous,
or equivalently, if one of the equivalent conditions (i) - (iv) in
Corollary \ruf{NE} is satisfied.

\Definition{} Two nil-polynomials $\6p\in\FF[W]$,
$\tilde{\6p}\in\FF[\tilde W]$ are called {\sl linearly (affinely)
equivalent} if there exists a linear (affine) isomorphism
$g:W\oplus\FF\to\tilde W\oplus\FF$ mapping $\Gamma_{\!\6p}$ onto
$\Gamma_{\!\tp}\,$.

\Proposition{HD} The nil-polynomials $\6p\in\FF[W]$,
$\tilde{\6p}\in\FF[\tilde W]$ are linearly equivalent if and only if
there exists a linear isomorphism $\alpha:W\to\tilde W$ and an
$\epsilon\in\GL(\FF)\cong\FF^{*}$ with
$\tilde\6p=\epsilon\circ\6p\circ\alpha^{-1}$.

\Proof Assume that $\6p$, $\tilde\6p$ are linearly equivalent.
Then there exist $\alpha\in\Hom(W,\tilde W)$,
$\beta\in\Hom(\FF,\tilde W)$ as well as $\gamma\in\Hom(W,\FF)$,
$\delta\in\FF$ such that $(x,t)\mapsto(\alpha x+\beta t,\gamma
x+\delta t)$ establishes a linear equivalence
$\Gamma_{\!\6p}\to\Gamma_{\!\tp}$. The ideal in $\FF(W\oplus\,\FF)$ of
all polynomials vanishing on $\Gamma_{\!\6p}$ is generated by
$t{-}\6p(x)$.  As a consequence we have for a suitable
$\epsilon\in\FF^{*}$
$$(\gamma x+\delta t)-\tilde\6p(\alpha x+\beta
t)=\epsilon(t-\6p(x))\Steil{for all}(x,t)\in W\oplus\FF\,.$$ Then
$\gamma=0$ and hence $\alpha$ is invertible.  Denote by
$\6q,\tilde\6q$ the quadratic parts of $\6p,\tilde\6p$. Then
$\tilde\6q(\alpha x+\beta t)=\epsilon\6q(x)$ for all $x,t$ implies
$\beta=0$ since the quadratic form $\tilde\6q$ is non-degenerate on
$\tilde W$.  Then $\tilde\6p(\alpha x)=\epsilon\6p(x)$ proves the
first claim. The converse is obvious.\qed

\noindent As a consequence, every equivalence class of nil-polynomials
in $\FF[W]$ is an orbit of the group $\FF^*\times\SL(W)$ acting in
the obvious way on $\FF[W]$. 

\medskip Corollaries \ruf{AC} and \ruf{NE} immediately imply
the following result.

\Proposition{UQ} Let $\6p$, $\tilde\6p$ be nil-polynomials associated
to the admissible algebras $\5N$, $\tilde\5N$. Then \0 $\5N$,
$\tilde\5N$ are isomorphic if and only if $\6p$, $\tilde\6p$ are
affinely equivalent. \1 In case $\5N$ has {\rm property (AH)} (for
instance, if $\5N$ admits a grading) then {\rm(i)} remains true with
`affinely' replaced by `linearly'.\Formend

\medskip For any pair of admissible algebras $\5N$, $\tilde{\5N}$ with
nil-polynomials $\6p\in\FF[W]$, $\tilde\6p\in\FF[\tilde W]$
Proposition \ruf{UQ} has the following obvious consequence.

\Corollary{PV} If $\5N$ has {\rm property (AH)} and $\5N$,
$\tilde{\5N}$ are isomorphic, then there exists an
$\epsilon\in\FF^{*}$ and a linear isomorphism $\alpha:W\to\tilde W$
with $\tilde\6p^{[k]}\circ\alpha=\epsilon\,\6p^{[k]}$ for all $k$, where
$\6p^{[k]}$ is the homogeneous part of degree $k$ in $\6p$.\Formend

\smallskip For $2\le k\le\nu(\5N)$ the homogeneous polynomial
$\6p^{[k]}$ is non-zero and gives a projective subvariety in the
projective space $\PP(\5N)$ associated to $\5N$. All these varieties
then are invariants for the algebra structure of $\5N$, provided $\5N$
has property (AH). It is worthwhile to mention that this remains true
for the leading homogeneous part also in the general situation, more
precisely:

\Proposition{VP} The statement of { \rm Corollary} \ruf{PV} remains
true for $k=\nu(\5N)$ even without requiring that $\5N$ has {\rm
property (AH)}.

\Proof Assume that $L:\5N\to\tilde{\5N}$ is an algebra
isomorphism. Then $\5N$, $\tilde{\5N}$ have the same nil-index, say
$\nu\ge2$. Without loss of generality we may assume that there are
pointings $\omega$, $\tilde\omega$ on $\5N$, $\tilde{\5N}$ with
$W=\ker(\omega)$, $\tilde W=\ker(\tilde\omega)$ and
$\6p=\omega\circ\exp_{2}$, $\tilde\6p=\tilde\omega\circ\exp_{2}$ on
$W,\tilde W$. Because of $L(\5N_{[1]})=\tilde{\5N}_{[1]}$ there is an
$\epsilon\in\FF^{*}$ with $\tilde\omega(L(a))=\epsilon\omega(a)$ for
every $a\in\5N_{[1]}$.  Further, there exists a linear isomorphism
$\alpha:W\to\tilde W$ and a linear map $\lambda:W\to\tilde{\5N}_{[1]}$
with $L(x)=\alpha(x)+\lambda(x)$ for all $x\in W$. But then
$\tilde\6p^{[\nu]}(\alpha x)=\tilde\omega((\alpha
x)^{(\nu)})=\tilde\omega((Lx)^{(\nu)})=
\tilde\omega(L(x^{(\nu)}))=\epsilon\omega(x^{(\nu)})=\epsilon\6p^{[\nu]}(x)$
since $x^{(\nu)}\in\5N_{[1]}$.  \qed

We illustrate by an example how \ruf{VP} can be applied to prove that
two given admissible algebras are not isomorphic: Anticipating
notation of Section \ruf{counter}, see \ruf{PQ}, consider
$\5M(Z^{3}+Y^{4}+X^{3}Z+X^{3}Y^{2}+X^{5}YZ)$ and
$\5M(Z^{3}+Y^{4}+X^{3}Z+X^{2}YZ+X^{3}Y^{2})$. These  are admissible
algebras of dimension 20 with nil-index 6, both having the same Hilbert
function $\{1,3,5,5,4,2,1\}$. It can be seen$^{2}$ that both algebras
do not have property (AH), compare also with Section
\ruf{counter}. Leading terms of nil-polynomials
$\6p,\tilde\6p\in\FF[x_{1},\dots,x_{20}]$ are, for instance,
$\6p^{[6]}=x_{1}^{4}x_{2}^{2}$ (for the first algebra) and
$\tilde\6p^{[6]}=x_{1}^{4}(19x_{1}^{2}-90x_{1}x_{2}+135x_{2}^{2})$
(for the second). Since the quadratic factor in $\tilde\6p^{[6]}$ is
not the square of a linear form we conclude with Proposition \ruf{VP}
that the algebras are not isomorphic.  

\medskip\noindent {\bf Remark.} There is a geometric interpretation of
the leading homogeneous term $\6p^{[\nu]}$: Identify $W$ in the
standard way with an affine open subset in the projective space
$\PP(\FF\oplus W)$. Hence $\PP(\FF\oplus W)=W\;\dot\cup\; \PP(W)$,
where $\PP(W)$ is the projective hyperplane at infinity.  The zero set
$T:=\{\6p=0\}\subset W$ is linearly equivalent to $\{\6f=0\}\cap
\varphi(W)$ where $\varphi:W\to \ker(\omega)\subset \5N$ is the linear
isomorphism from definition \ruf{RS}. Consider the Zariski closure
$\Cl(T)\subset \PP(\FF\oplus W)$. Then the set of points at infinity,
$T^\infty:=\Cl(T)\cap \PP(W)$, coincides with $\{[z]\in \PP(W):
\6p^{[\nu]}(z)=0\}$.  For a not algebraically closed field $\FF$ then
$T^\infty$ encodes in general less information then the homogeneous
part $\6p^{[\nu]}$. Indeed, for instance in case $\FF=\RR$ the
quadratic factor in $\tilde\6p^{[6]}$ above is positive definite on
$\RR^{2}$, that is, the zero locus of $\tilde\6p^{[6]}$ in the real
projective space $\PP_{19}(\RR)$ is the hyperplane $\{x_{1}=0\}$. This
suggests to consider projective varieties defined by the $\6p^{[\nu]}$
(or $\6p^{[k]}$, $k<\nu$) only in case of algebraically closed base
fields. In such a situation the corresponding divisor rather then the
mere zero set is an invariant equivalent to $\6p^{[\nu]}$.

\bigskip Every nil-polynomial $\6p$ associated to $\5N$ depends on
$\dim(\5N)-1$ variables. Another type of polynomial, closer to
\Ruf{DN}, can be defined as follows: 

\Definition{} The polynomial $\6f\in\FF[V]$ is called an {\sl extended
nil-polynomial} associated to $\5N$, if there exists a linear
isomorphism $\phi:V\to\5N$ and a pointing $\omega$ on $\5N$ such that
$\6f=\omega\circ\exp_{1}\!\circ\,\phi\,$.

It is clear that for the linear part $\6f^{[1]}\in\FF[V]$ of $\6f$ the
restriction of $\6f$ to $W:=(\6f^{[1]})^{-1}(0)$ is a nil-polynomial
associated to $\5N$ and that the graph $\Gamma_{\6p}\subset
W\oplus\FF$ is linearly equivalent to the hypersurface
$\6f^{-1}(0)\subset V$. Conversely, every nil-polynomial
$\6p\in\FF[W]$ associated to $\5N$ can be extended by
$\6f(x,t):=\6p(x)+t$ to an extended nil-polynomial
$\6f\in\FF(W\oplus\FF)$.

\bigskip For a given pointed algebra $\5(\5N,\omega)$ fix a
nil-polynomial $\6p=\omega\circ\exp_{2}\circ\,\phi\in\FF[W]$ in the
following and define the symmetric $k$-form $\omega_{k}$ on $W$ by
$$\omega_{k}(x_{1},x_{2},\dots,x_{k})=\omega\big((\phi x_{1})(\phi
x_{2})\cdots(\phi x_{k})\big)\,.\Leqno{HO}$$ Then we have the
expansion $\6p=\sum_{k\ge2}\6p^{[k]}$ into homogeneous parts, where
$$\6p^{[k]}(x)={1\over k!}\,\omega_{k}(x,\dots,x)\Leqno{NG}$$ and
$\omega_{2}$ is non-degenerate on $W$. Using $\6p^{[2]}$ and
$\6p^{[3]}$ we define a commutative (a priori not necessarily
associative) product $(x,y)\mapsto x\cd y$ on $W$ by
$$\omega_{2}(x\cd y,z)=\omega_{3}(x,y,z)\steil{for all}z\in
W\Leqno{HY}$$ and also a commutative product on $W\oplus\,\FF$ by
$$(x,s)(y,t):=(x\cd y,\omega_{2}(x,y))\,.\Leqno{HR}$$ For
$\5K:=\ker(\omega)$ there is a unique linear isomorphism
$\psi:\FF\to\5N_{[1]}$ such that $\pi=\psi\circ\omega$ is the
canonical projection $\5K\oplus\5N_{[1]}\to\5N_{[1]}$. With these
ingredients we have

\Proposition{HZ} With respect to the product \Ruf{HR} the linear map
$$W\oplus\,\FF\;\to\;\5N\,, \quad(x,s)\mapsto\phi(x)+\psi(s)\,,$$ is
an isomorphism of algebras. In particular, $W$ with product $x\cd y$
is isomorphic to the nilpotent algebra $\5N/\5N_{[1]}$ and has
nil-index $\nu(\5N){-}1$.

\Proof For all $x,y\in W$ we have

\smallskip\centerline{$(\phi(x)+\psi(s))(\phi(y)+\psi(t))=(N-A)+A$
with}

\smallskip\centerline{$N:=\phi(x)\phi(y)\in\5N\steil{and}A:=
\pi(\phi(x)\phi(y))=\psi(\omega_{2}(x,y))\in\5N_{[1]}$.}

\noindent It remains to show $N-A=\phi(x\cd y)$. But this follows from

\smallskip\centerline{$N-A\,\in\,\5K\Steil{and}\omega\big(\phi(x\cd
y)\phi(z)\big)=\omega\big(\phi(x)\phi(y)\phi(z)\big)=
\omega\big((N-A)\phi(z)\big)$} \noindent for all $z\in W$. \qed

\Corollary{HP} Every nil-polynomial $\6p$ on $W$ is uniquely
determined by its quadratic and cubic term, $\6p^{[2]}$ and
$\6p^{[3]}$. In fact, the other $\6p^{[k]}$ are recursively determined
by \Ruf{NG} and
$$\omega_{k+1}(x_{0},x_{1},\dots,x_{k})=\omega_{k}(x_{0}\cd
x_{1},x_{2},\dots,x_{k})\Leqno{HQ}$$ for all $k\ge2$ and
$x_{0},x_{1},\dots,x_{k}\in W$.\Formend

\medskip Corollary \ruf{HP} suggests the following question: Given a
 non-degenerate quadratic form $\6q$ and a cubic form $\6c$ on $W$.
 When does there exist a nil-polynomial $\6p\in\FF[W]$ with
 $\6p^{[2]}=\6q$ and $\6p^{[3]}=\6c\,$? Using $\6q$, $\6c$ we can
 define as above for $k=2,3$ the symmetric $k$-linear form
 $\omega_{k}$ on $W$ and with it the commutative product $x\cd y$ on
 $W$. A necessary and sufficient condition for a positive answer is
 that $W$ with this product is a nilpotent and associative algebra. As
 a consequence we get for every fixed non-degenerate quadratic form
 $\6q$ on $W$ the following structural information on the space of all
 nil-polynomials $\6p$ on $W$ with $\6p^{[2]}=\6q\,$: Denote by $\6C$
 the set of all cubic forms on $W$. Then $\6C$ is a linear space of
 dimension ${n+2\choose3}$, $n=\dim W$, and
$$\6C_{\6q}:=\{\6c\in\6C:\exists\steil{nil-polynomial $\6p$ on $W$
with}\6p^{[2]}=\6q,\6p^{[3]}=\6c\}\Leqno{QH}$$ is an algebraic
subset.

\KAP{Representations}{Representations of nil-algebras and adapted
decompositions}

In the following let $A$ be an arbitrary commutative nilpotent algebra
and $\2E$ a vector space. Also let $$\N:A\to\End(\2E)\;,\qquad
x\mapsto\N_{x}\;,$$ be an algebra homomorphism. For example, every
commutative algebra $A$ admits the faithful left-regular
representation $\6L:A\to \End(A^0)$, where $A^0$ is the unital
extension of $A$.  Consider the following characteristic subspaces of
$\2E$:
$$\2B:=\langle \N_{x}(\2E):x\in
A\rangle_{\FF}\Steil{and}\2K:=\bigcap_{x\in A}\ker(\N_{x})
\;.$$ Let us call every decomposition
$$\2E=\2E_{0}\oplus\2E_{1}\oplus\2E_{2}\oplus\2E_{3}\Steil{with}
\2B=\2E_{2}\oplus\2E_{3}\,,\;\2K=\2E_{0}\oplus\2E_{3}\Leqno{UR}$$
an {\sl \N-adapted decomposition} of $\2E$. It is obvious that
starting with $\2E_{3}:=\2B\cap\2K$ and choosing successively suitable
linear complements $\2E_2,$ $\2E_0$ and $\2E_1$ one always obtains
an \N-adapted decompositions for $\2E$.
Clearly, $\2E_{0}=0$ if the image algebra $\N(A)$
is maximal among all commutative nilpotent subalgebras of $\End(\2E)$.

\medskip Now assume that $\2E$ has finite dimension and that a
non-degenerate symmetric bilinear form $\6h:\2E\times\2E\to\FF$ is
fixed such that every $\N_{x}$ is selfadjoint with respect to $\6h\,$:
For every $T\in\End(\2E)$ the adjoint $T^{\star}$ is defined by
$\6h(Tv,w)=\6h(v,T^{*}w)$ for all $v,w\in\2E$. The orthogonal
`complement' of every linear subspace $\2L\subset\2E$ is
$\2L^{\perp}:=\{v\in\2E:\6h(v,\2L)=0\}$. The linear subspace $\2L$ is
called {\sl totally isotropic} if $\2L\subset\2L^{\perp}$.

\Proposition{HC} There exists an \N-adapted decomposition \Ruf{UR}
that is related to $\6h$ in the following way: \0 The three
subspaces $\2E_{0}$, $\2E_{1}\!\oplus\2E_{3}$ and $\2E_{2}$ are
mutually orthogonal with respect to $\6h$. \1 The subspaces $\2E_{1}$
and $\,\2E_{3}$ are totally isotropic, and hence have the same
dimension.

\Proof Choose an arbitrary $\N$-adapted decomposition
$\2E=\2E_{0}\oplus\tilde\2E_{1}\oplus\2E_{2}\oplus\2E_{3}$. Since
every $\N_{x}$ is selfadjoint we have $\2K=\2B^{\perp}$. In
particular, $\2E_{3}$ is totally isotropic. We get further
$\dim(\2E)=\dim(\2K)+\dim(\2B)\,$,
$\dim(\tilde\2E_{1})=\dim(\2E_{3})$,
$\,\2E_{0}\oplus\2E_{2}\oplus\2E_{3}\; =\;\2E_{3}^{\perp}$ and that
the three spaces $\2E_{0}$, $\2E_{2}$, $\2E_{3}$ are mutually
orthogonal.\nline Further we conclude from
$\2E_{2}\!\cap\2E_{2}^{\perp}\,\subset\;\2B^{\perp}=\2K\,$ that
$\2E_{2}\!\cap\2E_{2}^{\perp}\,\subset\;\2E_{2}\cap\2K=0$. In the same
way we conclude from
$\2E_{0}\!\cap\2E_{0}^{\perp}\,\subset\;\2K^{\perp}=\2B\,$ that
$\2E_{0}\!\cap\2E_{0}^{\perp}\,\subset\;\2E_{0}\!\cap\2B=0$.  As a
consequence we get
$\;(\2E_{0}^{\perp}\!\cap\2E_{2}^{\perp})\cap\,\2E_{3}^{\perp}=\2E_{3}$
and $\2E_{0}^{\perp}+\2E_{2}^{\perp}=\2E\,$. Now choose a linear
subspace $\2E_{1}\subset(\2E_{0}^{\perp}\cap\2E_{2}^{\perp})$ with
$\2E_{1}\oplus\2E_{3}=\2E_{0}^{\perp}\!\cap\2E_{2}^{\perp}$. Counting
dimensions we get $\dim(\2E_{1})=\dim(\tilde\2E_{1})$ from
$\dim(\2E_{0}^{\perp}\cap\2E_{2}^{\perp})=\dim(\2E_{0}^{\perp})+
\dim(\2E_{2}^{\perp})-\dim(\2E_{0}^{\perp}+\2E_{2}^{\perp})$ and thus
that $\2E=\2E_{0}\oplus\2E_{1}\oplus\2E_{2}\oplus\2E_{3}$ is an
$\N$-adapted decomposition satisfying (i). The form $\6h$ is
non-degenerate on $\2E_{1}\oplus\2E_{3}$. Because of
$\dim(\2E_{1})=\dim(\2E_{3})$ we finally may assume without loss of
generality that also $\2E_{1}$ is totally isotropic.\qed

We call every \N-adapted decomposition satisfying (i), (ii) above an
{\sl $(\N,\6h)$-adapted decomposition} of the representation space
$\2E$.  In the following we give two applications: 

\noindent Let $\5N$ be an admissible algebra with pointing
$\omega$. Clearly, in general the quotient $B:=\5N/\5N_{[1]}$ is a
non-admissible nilpotent algebra, say with product $(x,y)\mapsto
x\bullet y$.  Left multiplication yields a (non-faithful)
representation $\N:B\to\End(B)$ in terms of the multiplication
operator $\N_{x}:y\mapsto x\bullet y$.  Further, the symmetric
bilinear form $\6b(x,y)=\omega(xy)$ on $\5N$ factors to a
non-degenerate symmetric bilinear form $\6h$ on $\5N/\5N_{[1]}$ and
all $\N_{x}$ are selfadjoint with respect to $\6h$.  Note that $\5N$
is isomorphic to $(\5N/\5N_{[1]})\times \FF$ with multiplication given
by
$$ (x,a)\diamond
(y,b):=(\N_x(y),\6h(x,y))=(\N_y(x),\6h(y,x))\,.\Leqno{JG}$$ Instead of
$\5N/\5N_{[1]}$ we use the isomorphic algebra $W:=\ker(\omega)$ with
product $x\cd y$, as given in Proposition \ruf{HZ}. Then the form
$\6h$ is the restriction of $\6b$ to $W$. If $\5N$ has nil-index
$\nu\ge2$ then the subalgebra $\N(W)\subset\End(W)$ has nil-index
$\nu{-}2$. If $W=W_{0}\oplus W_{1}\oplus W_{2}\oplus W_{3}$ is a
$(W,\6h)$-adapted decomposition, then $\5N':=W_{1}\oplus W_{2}\oplus
W_{3}\oplus\5N_{[1]}$ and $\5N'':=W_{0}\oplus\5N_{[1]}$ are admissible
subalgebras with $\nu(\5N'')\le2\,$, and $\5N$ is a smash product of
$\5N'$ with $\5N''$, as defined in Section $\ruf{examples}$.

\Proposition{HI} Every admissible algebra of nil-index $\le3$ has
a grading.

\Proof As indicated above $\5N$ is isomorphic to $W\times \FF$ with
product \ruf{JG} where $W=\ker(\omega)\subset\5N$ is the nilpotent
subalgebra isomorphic to $\5N/\5N_{[1]}$.  Let $\6N:W\to \End(W)$ as
above and consider the subalgebra $\N(W)\subset\End(W)$. Let
$W=W_{0}\oplus W_{1}\oplus W_{2}\oplus W_{3}$ be a $(W,\6h)$-adapted
decomposition. Since $\N(W)$ has nil-index $\le1$ we have
$W_{2}=0$. Put $\5N_{2}:=W_{1}$, $\,\5N_{4}:=W_{3}$,
$\,\5N_{6}:=\Ann(\5N)$ and $\5N_{3}:=W_{0}$. Then
$\5N=\5N_{2}\oplus\5N_{3}\oplus\5N_{4}\oplus\5N_{6}$ is a grading of
$\5N$.\qed

The estimate $\nu(\5N)\le3$ in Proposition \ruf{HI} is sharp as a
counterexample in Section \ruf{counter} with nil-index 4 and dimension
8 will show.

\medskip The next result improves Proposition \ruf{ZT} in the case of
admissible algebras.

\Proposition{QI} For every admissible algebra $\5N$ and every
$\pi\in\Pi(\5N)$ there exists a subgroup of $\Aff(S_{\pi})$ acting
transitively on $S_{\pi}\cap\5N_{[4]}$. In particular, $\5N$ has {\rm
Property (AH)} if $\5N$ has nil-index $\le4$.

\Proof Put $S:=S_{\pi}$ as shorthand and denote by $\6h$ the
restriction of $\6b_{\pi}$ to $W:=\ker(\pi)$. As above,
$\N_{x}\in\End(W)$ is the multiplication operator $y\mapsto x\cd
y$. Choose a $(W,\6h)$-adapted decomposition $W=W_{0}\oplus
W_{1}\oplus W_{2}\oplus W_{3}$ and denote by $\pi_{k}\in\End(W)$ the
canonical projection with range $W_{k}$ for $0\le k\le3$. Then
$(W_{0}+W_{3})\cd W=0$, $W\cd W\subset W_{2}\oplus W_{3}$ and $W'\cd
W_{2}\subset W_{3}$, where $W':=\5N_{[4]}\cap W$.

Now fix a point $a\in S\cap W'$. Because of Proposition \ruf{ZT} it is
enough to show $g(a)\in\5N_{[3]}$ for some
$g\in\Aff(S)\cap\Aff(\5N_{[4]})\,$: Put $\6P:=\pi_{3}\circ
\N_{c}\circ\pi_{2}-\pi_{2}\circ \N_{c}\circ\pi_{1}$ for
$c:=a-\pi(a)\in W'$. Then $\6P^{*}=-\6P$ and
$\6Q:={1\over2}\N_{c}+{1\over6}\6P\in\End(W)$ is nilpotent with
$\6Q(W)\subset\5N_{[3]}$. Set $\5N^{\,0}:=\FF\One\oplus
W\oplus\5N_{[1]}$ and define $\lambda\in\End(\5N^{\,0})$ by
$$(s\One,x,t)\mapsto(0,\6Qx-sc,\6h(c,x))\steil{for all}s\in\FF,x\in
W,\,t\in\5N_{[1]}\,.$$ $\lambda$ is nilpotent and maps $\5N^{\,0}$ to
$\5N_{[4]}$. Therefore the unipotent operator
$g:=\exp(\lambda)\in\GL(\5N^{\,0})$ exists. Clearly, $\5U:=\One+\5N$
and $\One+\5N_{[4]}$ are $g$-invariant.

We proceed as in the proof of Theorem \ruf{NB}. Denote by $\5F$ the
space of all polynomial maps $\5N^{\,0}\to\5N_{[1]}$ of degree $\le 4$
and define the nilpotent operator $\xi\in\End(\5F)$ by
$$(\xi f)(z):=f'(z)(\lambda z)\Steil{for
all}f\in\5F\steil{and}z\in\5N^{\,0}\;.$$ We claim that $g$
leaves $\One+S$ invariant. For this we only have to show that
$\xi\hat\6f$ vanishes on $\One+S$, where $\hat\6f\in\5F$ is defined by
$\hat\6f(s\One,x,t):=t+\pi(\exp_{2}x)$.  But this just means that
$$\6h(c-\6Qx,x+{\scriptstyle{1\over2}}\6N_{x}x+
{\scriptstyle{1\over6}}\6N_{x}^{2}x)=\6h(c,x)$$ holds for all $x\in
W$, or equivalently
$$\6h(2\6Qx-\6N_{c}x,x)\;=\;\6h(3\6Qx-\6N_{c}x,\6N_{x}x)\;=\;
\6h(\6Qx,\6N_{x}^{2}x)\;=0\;\,.$$ The first term vanishes since
$2\6Q-\6N_{c}$ is skew adjoint. The two other terms vanish since
${3\6Qx-\6N_{c}x}\in\5N_{[2]}$ and $\6Qx\in\5N_{[3]}$.
This proves the claim, and as a consequence we
get $g(\One+a)=\One+b\,$ with $b\in (c\cd W+
W_{3}+\5N_{[1]})\subset\5N_{[3]}$.\qed

\medskip We conclude the section with a lower bound for the dimension
of $\der(\5N)$ and $\Aut(\5N)$.  For fixed $\pi\in\Pi(\5N)$, by the
proof of Lemma \ruf{PA} every $\lambda\in\Hom(\5N,\5N_{[1]})$ with
$\lambda(\5N_{[1]})=0$ is of the form $\lambda=\pi\circ M_{b}$,
$b\in\ker\pi$. Therefore $c\mapsto\pi\circ M_{c}$ induces for every
$k>0$ a linear isomorphism
$\5N_{[k]}/\5N_{[1]}\;\cong\;\Hom(\5N/\5N^{k},\5N_{[1]})$. This
implies$$\dim\5N_{[k]}+\dim\5N^{k}=\dim\5N+1\Steil{for all}0\le
k\le\nu(\5N)\,.$$ As a consequence we get from Proposition \ruf{QI}
the

\Corollary{} $\dim\Aut(\5N)\;\ge\;\dim\5N/\5N^{4}$ for every
admissible algebra $\5N$. The same lower bound holds for the dimension
of $\,\der(\5N)$.

\Proof In case $\nu(\5N)<4$ the algebra $\5N$ is gradable by
Proposition \ruf{HI} proving the claim.  In case
$\nu(\5N)\ge4$ Proposition \ruf{QI} together with
$\dim\,(\5N_{[4]}/\5N_{[1]})=\dim(\5N/\5N^{4})$ implies
$$\dim\Aff(S_{\pi})\;\ge\;\dim\5N/\5N^{4}\,+\,\dim\GL(S_{\pi})\,,
\;\;\rlap{ that is}$$ 
$$\dim\Aut(\5N)\;\ge\;\dim\5N/\5N^{4}+
\dim\,\{g\in\Aut(\5N):g\circ\pi=\pi\circ g\}$$ as a consequence of
Proposition \ruf{NI}. The estimate for $\der(\5N)$ follows with
\ruf{UO}.\qed

\KAP{examples}{Some examples}

For every admissible algebra $\5N$ and associated nil-polynomial $\6p$
the degree of $\6p$ is the nil-index of $\5N$, provided
$\dim\5N>1$. In the following we give a method how large classes of
nil-polynomials of degrees 3 and 4 can be constructed.

Nil-polynomials of degree 2 on a given vector space $W\ne0$ are quite
obvious -- these are precisely all non-degenerate quadratic forms
$\6q$ on $W$. Then $\6q$ is associated to the admissible algebra
$\5N:=W\oplus\FF$ whose commutative product is uniquely determined by
$(x,t)^{2}=(0,\6q(x))$ for all $x\in W$ and $t\in\FF$ (special case of
\ruf{JG}). Clearly, $\5N$ is canonically graded by $\5N_{1}=W$ and
$\5N_{2}=\FF$. Notice that in every dimension $n\ge2$ the number of
isomorphy classes for admissible algebras with nil-index 2 is: 1 if
$\FF$ is algebraically closed; $[n/2]$ if $\FF$ is the real field;
infinite if $\FF$ is the rational field.

For every pair of nil-polyno\-mials $\6p'\in\FF[W']$,
$\6p''\in\FF[W'']$ we get a new nil-polynomial
$\6p:=\6p'\oplus\6p''\in\FF[W'\oplus W'']$ by setting
$\6p(x',x''):=\6p'(x')+\6p''(x'')$ for all $x'\in W'$ and $x''\in
W''$. Let us call a nil-polynomial {\sl reduced} if it is not affinely
equivalent to a direct sum $\6p'\oplus\6p''$ of nil-polynomials with
$\6p''$ of degree two. Then it is clear that every nil-polynomial
$\6p$ is affinely equivalent to a direct sum $\6p'\oplus\6p''$ with
$\6p'$ reduced and $\6p''$ of degree $\le2$. The nil-polynomials
$\6p'$, $\6p''$ are uniquely determined up to affine equivalence.

On the level of algebras the direct sum of nil-polynomials corresponds
to the following construction: Let $(\5N',\omega')$,
$(\5N'',\omega'')$ be pointed algebras. Then
$\5I:=\{(s,t)\in\5N_{[1]}'\times\5N_{[1]}'':\omega'(s)+\omega''(t)=0\}$
is an ideal in $\5N'\times\5N''$ and
$\5N:=(\5N'\times\5N'')/\5I$ becomes a pointed algebra with
respect to the pointing
$(s,t)+\5I\;\mapsto\;\omega'(s)+\omega''(t)$. We write
$\5N=\5N'\vee\5N''$ and call it the {\sl smash product} of the
pointed algebras $\5N_{j}$. This product is commutative and
associative in the sense that $(\5N_{1}\vee\5N_{2})\vee\5N_{3}$ and
$\5N_{1}\vee(\5N_{2}\vee\5N_{3})$ are canonically isomorphic. We call
the admissible algebra $\5N$ {\sl reduced} if it is not isomorphic to
a smash product $\5N_{1}\vee\5N_{2}$ with $\nu(\5N_{2})=2$. Notice
that every admissible algebra $\5N'$ of dimension one is neutral with
respect to the smash product, that is, $\5N\vee\5N'\cong\5N$. Notice
also that any smash product of gradable admissible algebras is also
gradable.

\medskip
\centerline{\bf Nil-polynomials of degree 3}

Let $p\in\FF[V]$ be a homogeneous polynomial of degree $d\ge2$, that
is, $p(x)=q(x,\dots,x)$ for all $x\in V$ and a uniquely determined
symmetric $d$-linear form $q:V^{d}\to\FF$. Then $p$ is called {\sl
non-degenerate} if $V_{p}^{0}=0\;$ for $$V_{p}^{0}:=\{a\in
V:q(a,V,\dots,V)=0\}\;=\;\{a\in V:p(x+a)=p(x)\steil{for all}x\in
V\}\;.$$

\Proposition{HU} Let $W$ be an $\FF$-vector space of finite dimension
and $\6q$ a non-degenerate quadratic form on $W$. Suppose furthermore
that $W=W_{1}\oplus W_{2}$ for totally isotropic (with respect to
$\6q$) linear subspaces $W_{k}$ and that $\6c$ is a cubic form on
$W_{1}$. Then, if we extend $\6c$ to $W$ by $\6c(x+y)=\6c(x)$ for all
$x\in W_{1}$, $y\in W_{2}$, the sum $\6p:=\6q+\6c$ is a nil-polynomial
on $W$, and $\6p$ is reduced if and only if $\6c$ is non-degenerate on
$W_{1}$. For all cubic forms $\6c$, $\tilde\6c$ on $W_{1}$ with
nil-polynomials $\6p=\6q+\6c$, $\;\tilde\6p=\6q+\tilde\6c$ and
associated admissible algebras $\5N$, $\tilde\5N$ the following
conditions are equivalent:\0 $\tilde\6c=\6c\circ g$ for some
$g\in\GL(W_{1})$.\1 $\5N$, $\tilde\5N$ are isomorphic as algebras.

\Proof $\omega_{3}(x,y,t)=0$ for all $t\in W_{2}$ implies $W\!\cd
W\subset W_{2}$ and $W\!\cd W_{2}=0$, that is, $(x\cd y)\cd z=0$ for
all $x,y,z\in W$, see \Ruf{HY}. This implies that $\6p:=\6q+\6c$ is a
nil-polynomial. Fix a decomposition $W_{1}=W_{1}'\oplus W_{1}''$ with
$W_{1}'':=\{a\in W_{1}:\6c(x+a)=\6c(x)\steil{for all}x\in
W_{1}\}$. Then also $W_{1}''=\{a\in W_{1}:a\cd W_{1}=0\}$ and we put
$$W_{2}':=\{y\in W_{2}:y\perp W_{1}''\}\Steil{and} W_{2}'':=\{y\in
W_{2}:y\perp W_{1}'\}\;.$$ Then for $W':=W_{1}'\oplus W_{2}'$ and
$W'':=W_{1}''\oplus W_{2}''$ we have the orthogonal decomposition
$W=W'\oplus W''$ with $W\!\cd W''=0$ and $W\!\cd W\subset
W_{2}'$. Denote by $\6p'$ and $\6p''$ the restriction of $\6p$ to $W'$
and $W''$ respectively. Then $\6p=\6p'\oplus\6p''$ as direct sum of
nil-polynomials and $\6p''$ has degree $\le2$.  These considerations
show that $\6p$ is reduced if and only if $W_{1}''=0$. It remains to
verify the equivalence of (i) and (ii).

\noindent \To12 There exists a unique $g^{\sharp}\in\GL(W_{2})$ with
$\omega_{2}(gx,y)=\omega_{2}(x,g^{\sharp}y)$ for all $x\in W_{1}$ and
$y\in W_{2}$. But then $\tilde\6p=\6p\,\circ h$ for $h:=g\times
(g^{\sharp})^{-1}\in\O(\6q)\subset\GL(W)$. By Proposition \ruf{UQ} the
algebras $\5N$, $\tilde\5N$ are isomorphic.\nline \To21 By Proposition
\ruf{HU} the algebra $\5N/\5N_{[1]}$ is isomorphic to $W$ with the
product $x\cd y$ determined by $\6p^{[2]}=\6q$ and
$\6p^{[3]}=\6c$. Choose, as above, decompositions $W=W'\oplus W''$,
$W'=W_{1}'\oplus W_{2}'$ and $W''=W_{1}''\oplus W_{2}''$ . Then the
cubic form $\6c$ `essentially lives' on the subspace $W_{1}'\subset
W_{1}$ and $W_{1}''\oplus W_{2}$ is the annihilator of $W$. In the
same way the product $x\tpp y$ determined by $\tilde\6p^{[2]}=\6q$ and
$\tilde\6p^{[3]}=\tilde\6c$ on $W$ gives an algebra $\tilde W$
isomorphic to $\tilde\5N/\tilde\5N_{[1]}$. We choose again
decompositions $\tilde W=\tilde W'\oplus\tilde W''$, $\,\tilde
W'=\tilde W_{1}'\oplus\tilde W_{2}'$ and $\tilde W''=\tilde
W_{1}''\oplus\tilde W_{2}''$ and have $\tilde W_{1}''\oplus\tilde
W_{2}=\Ann(\tilde W)$. \nline Now assume (ii) and fix an algebra
isomorphism $h:W\to\tilde W$. Then there exists a linear isomorphism
$\alpha:W_{1}'\to\tilde W_{1}'$ with $\;h(x)\equiv
\alpha(x)\,\mod\,\Ann(W)\;$ for all $x\in W_{1}$. Let
$\beta:W_{1}''\to\tilde W_{1}''$ be an arbitrary linear isomorphism
and consider $g:=\alpha\oplus\beta$ as an element of
$\GL(W_{1})$. Replacing $\tilde\6c\,$ by $\,\tilde\6c\circ g^{-1}$ we
may assume without loss of generality that $g=\id_{W_{1}}$. But then
$x\cd x=x\tpp x$ and thus $6\,\6c(x)=\omega_{2}(x\cd
x,x)=\omega_{2}(x\tpp x,x)=6\,\tilde\6c(x)$ for all $x\in W_{1}$.\qed

\medskip Suppose that $W_{1}\cong\FF^{m}$ with coordinates
$x=(x_{1},\dots,x_{m})$ has dimension $m>0$ in Proposition
\ruf{HU}. Then $W\cong\FF^{2m}$ with coordinates
$(x,y)=(x_{1},\dots,x_{m},y_{1},\dots,y_{m})$ and we may assume
$\6q(x,y)=x_{1}y_{1}+\dots+x_{m}y_{m}$. As already mentioned, the
linear space $\6C$ of all cubic forms on $W_{1}$ has dimension
${m+2\choose3}$. The group $\GL(W_{1})$ acts on $\6C$ from the right
and has dimension $m^{2}$ over $\FF$. The difference of dimensions is
$m\choose3$.  But this number is also the cardinality of the subset
$J\subset\NN^{3}$, consisting of all triples $j=(j_{1},j_{2},j_{3})$
with $1\le j_{1}<j_{2}<j_{3}\le m$. Consider the affine map
$$\alpha:\FF^{J}\to\6C\,,\quad (t_{j})\mapsto \6c^{}_{0}+\sum_{j\in
J}t_{j}\6c_{j}\,,\Leqno{HX}$$ where
$\6c_{0}:=x_{1}^{3}+\dots+x_{m}^{3}$ and
$\6c_{j}:=x_{j_{1}}x_{j_{2}}x_{j_{3}}$ for all $j\in J$. Notice that
the nil-polynomial
$\6q+\6c_{0}=\bigoplus_{k=1}^{m}(x_{k}y_{k}+x_{k}^{3})$ is a direct
sum of $m$ nil-polynomials, that is, the corresponding admissible
algebra is an $m$-fold smash power of the cyclic algebra $\FF
t\oplus\FF t^{2}\oplus\FF t^{3}$, where $t^{4}=0$.

In case $\FF=\RR$ or $\FF=\CC$, for a suitable neighbourhood $U$ of
$0\in\FF^{J}$ the map $\alpha:U\to\6C$ intersects all
$\GL(n,\FF)$-orbits in $\6C$ transversally. Indeed, since all
partial derivatives of $\6c_{0}$ are monomials containing a square,
the tangent space at $\6c_{0}$ of its $\GL(n,\FF)$-orbit is
transversal to the linear subspace $\langle\6c_{j}:j\in
J\rangle_{\FF}$ of $\6C$. In particular, in case $m\ge3$ there is a
family of dimension ${m\choose3}\ge1$ over $\FF$ $(=\RR\hbox{ or
}\CC)$ of pairwise different $\GL(n,\FF)$-orbits and thus of
non-equivalent nil-polynomials of degree $3$ on $W$. Notice that in
case $m=3$ the mapping $\alpha$ in \Ruf{HX} reduces to
$$\alpha:\FF\to\6C,\quad t\mapsto
x_{1}^{3}+x_{2}^{3}+x_{3}^{3}+tx_{1}x_{2}x_{3}\,.\Leqno{YH}$$

\Corollary{} Let $\FF$ be either $\RR$ or $\CC$. Then in every
dimension $n\ge7$ there is an infinite (in fact uncountable) number of
isomorphy classes of admissible algebras of dimension $n$ and
nil-index 3.\Formend

Proposition \ruf{HC} together with Proposition \ruf{HU} generalizes
Theorems 3.3 and 4.1 in \Lit{ELRO} from the case of algebraically
closed base fields to arbitrary fields of characteristic 0, see
Proposition \ruf{UT} below. For every admissible algebra $\5N$ with
nil-index $\nu$ the non-degeneracy of $\6b_{\pi}$ on $\5N/\5N^{\nu}$
implies $\dim(\5N^{k}/\5N^{\nu})\le\dim(\5N/\5N^{j})$ for all
$\nu/2<k<\nu$ and $j=\nu-k+1$ and thus $H(\nu-1)\le H(1)$ for the
Hilbert function $H=H_{\5N}$. This means in the special case of
nil-index $\nu=3$ that the Hilbert function of $\5N$ has the form
$\{1,m,n,1\}$ with $m\ge n\ge1$.

\Proposition{UT} {\bf (Classification of admissible algebras with
nil-index 3)} The admissible algebras $\5N$ of nil-index 3 are, up to
isomorphism, precisely the smash products $\5N'\vee\5N''$, where
$\5N'$ is a reduced algebra of dimension $>1$ as described in
Proposition \ruf{HU} and $\5N''$ has nil-index at most two. If $\5N$
has Hilbert function $\{1,m,n,1\}$ then $\5N'$ has symmetric Hilbert
function $\{1,n,n,1\}$ and admits a canonical grading. Furthermore,
the Hilbert function of $\5N''$ is $\{1,m-n,1\}$ if $m>n$, and is
$\{1,1\}$ if $m=n$ (that is, $\dim(\5N'')=1$, or equivalently,
$\5N=\5N'$ in this case).\Formend

Notice that every admissible algebra with a canonical grading has
symmetric Hilbert function \Ruf{UH}. The converse is not true, see e.g.
the example of dimension 23 in Section \ruf{counter}.

\bigskip\noindent{\bf Moduli algebras of type $\bf\tilde\3E_{6}\,$.}
~With Proposition \ruf{HZ} it is possible to compute for a given
nil-polynomial an admissible algebra it is associated to. For the
nil-polynomials in Proposition \ruf{HU} in case $\dim W=6$ there is a
connection to moduli algebras associated to simple elliptic
singularities $\tilde E_{6}$, see \Lit{EAWO}, \p306: Let $\FF=\CC$ and
for $t\in\CC$ with $t^{3}+27\ne0$ consider the nilpotent algebra
$\5N_{t}=\5M(X^{3}+Y^{3}+Z^{3}+tXYZ)$, compare the notation \Ruf{PQ}
below. Then, with $x,y,z$ being the residue classes of $X,Y,Z$ a basis
for $\5N_{t}$ is $x,y,z,yz,xz,xy,xyz$ with $\Ann(\5N_{t})=\CC xyz$ and
$$\6p_{t}=tx_{1}^{3}+tx_{2}^{3}+tx_{3}^{3}-18x_{1}x_{2}x_{3}+
x_{1}x_{4}+x_{2}x_{5}+x_{3}x_{6}\Leqno{GH}$$ is a nil-polynomial
associated to $\5N_{t}$. Notice that the cubic part of $\6p_{t}$
occurs already in (3.1) of \Lit{EAWO}. Notice also that $W_{1}=\langle
x,y,z\rangle_{\CC}$ and $W_{3}=\langle yz,xz,xy\rangle_{\CC}$ are two
totally isotropic subspaces as occurring in Proposition
\ruf{HU}. Although for the moduli algebras the three parameters $t$
with $t^{3}+27=0$ have to be excluded (for these three values of $t$
the singularity of $\{X^3+Y^3+Z^3+tXYZ=0\}$ is not isolated),
$\6p_{t}$ is a nil-polynomial also for these $t$ and is associated to
an admissible algebra of nil-index 3 as well. For $t\ne0$ and
$s:=-18/t$ it is easy to see that the nil-polynomial $\6p_{t}$ in
\Ruf{GH} is linearly equivalent to
$$x_{1}^{3}+x_{2}^{3}+x_{3}^{3}+sx_{1}x_{2}x_{3}+
x_{1}x_{4}+x_{2}x_{5}+x_{3}x_{6}\,,\hbox{ compare with \Ruf{YH}}\,.$$

\bigskip 
\medskip \centerline{\bf Nil-polynomials of degree 4}

The method in Proposition \ruf{HU} can be generalized to get
nil-polynomials of higher degrees, say of degree 4 for simplicity.
Throughout the subsection we use the notation \Ruf{KP}. Because of
Propositions \ruf{UQ} and \ruf{QI} it is not necessary to distinguish
between linear and affine equivalence for nil-polynomials of degree
4.

For fixed $n,m\ge1$ let $W=W_{1}\oplus W_{2}\oplus W_{3}$ be a vector
space with $W_{1}=\FF^{n}$, $W_{2}=\FF^{m}$ and let $\6q$ be a fixed
non-degenerate quadratic form on $W$ in the following. Assume that
$W_{1}$, $W_{3}$ are totally isotropic and that $W_{1}\oplus W_{3}$,
$W_{2}$ are orthogonal with respect to $\6q$. Then $W$ has dimension
$2n+m$, and without loss of generality we assume that for suitable
$\epsilon_{1},\dots,\epsilon_{m}\in\FF^{*}$
$$\6q(y)=\sum_{k=1}^{m}\epsilon_{k}y_{k}^{(2)}\Steil{if}y\in W_{2}\,.$$ As
before let $\6C$ be the space of all cubic forms on $W$. Our aim is to
find cubic forms $\6c\in\6C_{\6q}$ that are the cubic part of a
nil-polynomial of degree 4.

Denote by $\C$ the space of all cubic forms $\6c$ on $W_{1}\oplus
W_{2}$ such that $\6c(x+y)$ is quadratic in $x\in W_{1}$ and linear in
$y\in W_{2}$, or equivalently, which are of the form
$$\6c(x+y)={1\over2}\sum_{k=1}^{m}\sum_{i,j=1}^{n}\6c_{ijk}x_{i}x_{j}y_{k}
\Steil{for all}x\in W_{1},y\in W_{2}$$ with suitable coefficients
$\6c_{ijk}=\6c_{jik}\in\FF$.  Extending every $\6c\in\C$ trivially to
a cubic form on $W$ we consider $\C$ as a subset of $\6C$.

For fixed $\6c\in\C$ the symmetric 2- and 3-linear forms
$\omega_{2},\omega_{3}$ on $W$ are defined by
$\omega_{2}(x,x)=2\6q(x)$ and $\omega_{3}(x,x,x)=6\6c(x)$ for all
$x\in W$. With the commutative product $x\cd y$ on $W$, see \Ruf{HY},
define in addition also the $k$-linear forms $\omega_{k}$ by \Ruf{HQ}
for all $k\ge4$. Then,
for every $x,y\in W_{1}$ the identity
$\omega_{2}(x\cd y,t)=\omega_{3}(x,y,t)=0$ for all $t\in W_{1}\oplus
W_{3}$ implies $x\cd y\in W_{2}$, that is $W_{1}\cd W_{1}\subset
W_{2}$. In the same way $\omega_{2}(x\cd y,t)=0$ for all $x\in W_{1}$,
$y\in W_{2}$ and $t\in W_{2}\oplus W_{3}$ implies $W_{1}\cd
W_{2}\subset W_{3}$. Also $W_{j}\cd W_{k}=0$ follows for all $j,k$
with $j+k\ge4$. Therefore $\6c$ belongs to $\6C_{\6q}$ if and only if
$(a\cd b)\cd c$ is symmetric in $a,c\in W_{1}$ for every $b\in W_{1}$.

In terms of the standard basis $e_{1},\dots,e_{m}$ of $W_{2}=\FF^{m}$
we have
$$a\cd b=\sum_{k=1}^{m}\Big(\sum_{i,j=1}^{n}\epsilon_{k}^{-1}
\6c_{ijk}a_{i}b_{j}\Big) e_{k}\Steil{for all}a,b\in W_{1}$$ and thus
with $\Theta_{i,j,r,s}:=\sum_{k=1}^{m}\epsilon_{k}^{-1}
\6c_{ijk}\6c_{rsk}$ we get the identity
$$\omega_{2}\big((a\cd b)\cd
c,d\big)=\sum_{i,j,r,s=1}^{n}\!\!\Theta_{i,j,r,s}\,a_{i}
b_{j}c_{r}d_{s}\Steil{for all}a,b,c,d \in W_{1}\,,\hbox{ implying}$$
\vskip-15pt$$\6A:=\C\cap\6C_{\6q}=\{\6c\in\C:\Theta_{i,j,r,s}\steil{is
symmetric in}i,r\}\,.\Leqno{UF}$$ Notice that the condition in
\Ruf{UF} implies that $\Theta_{i,j,r,s}$ is symmetric in all
indices. $\6A$ is a rational subvariety of the linear space $\C$, it
consists of all those $\6c$ for which the corresponding product $x\cd
y$ on $W$ is associative. The group
$\Gamma:=\GL(W_{1})\times\O(\6q_{|W_{2}})\;\subset\;\GL(W_{1}\oplus
W_{2})$ acts on $\C$ by $\6c\mapsto\6c\circ\gamma^{-1}$ for every
$\gamma\in\Gamma$. Furthermore, $(g,h)\mapsto(g,h,(g^{\sharp})^{-1})$
embeds $\Gamma$ into $\O(\6q)$, compare the proof of Proposition
\ruf{HU}. As a consequence, the subvariety $\6A\subset\C$ is invariant
under $\Gamma$.

\medskip We are only interested in the case where $\6c\in\6A$ is
non-degenerate on $W_{1}\oplus W_{2}$. Then $W_{k+1}=\langle W_{1}\cd
W_{k}\rangle_{\FF}$ holds for $k=1,2$, implying
$m\le{n+1\choose2}$. Put $\5N:=\bigoplus_{k>0} W_{k}$ with
$W_{4}:=\FF$ and $W_{k}:=0$ for all $k>4$. The product \Ruf{HR}
realizes $\5N$ as graded admissible algebra. The non-degeneracy of
$\6c$ gives in addition $\5N^{k}=\bigoplus_{\ell\ge k}\!W_{\ell}$ for
all $k>0\,$, that is, the grading is canonical.  Conversely, every
canonically graded admissible algebra of nil-index 4 occurs (up to
isomorphism) this way with a non-degenerate $\6c$ as above. Further
admissible algebras with nil-index 4 can be obtained from $\5N$ as
above by taking $\5N\vee\5N'$ with $\5N'$ an arbitrary admissible
algebra of nil-index $\le3$. But all these algebras are gradable by
Proposition \ruf{HI}. Since there exist non-gradable admissible
algebras of nil-index 4 (see Section \ruf{counter}) the classification
problem in the nil-index 4 case must be more involved than the one in
Proposition \ruf{UT}.

\smallskip Let us consider the special case $n=2$ with
$m={n+1\choose2}=3$ in more detail (among these are in case $\FF=\CC$
also all nil-polynomials of maximal ideals of moduli algebras
associated to singularities of type $\tilde E_{7}$, see \Lit{EAWO},
\p307). For simplicity we assume that for suitable coordinates
$(x_{1},x_{2})$ of $W_{1}$, $(y_{1},y_{2},y_{3})$ of $W_{2}$ and
$(z_{1},z_{2})$ of $W_{3}$ the quadratic form $\6q$ is given by
$$\6q=x_{1}z_{1}+x_{2}z_{2}+y_{1}^{(2)}+ y_{2}^{(2)}+\epsilon
y_{3}^{(2)}\Steil{for fixed}\epsilon\in\FF^{*}\Leqno{KU}$$ (in case
$\FF=\RR,\CC$ this is not a real restriction). For every $t\in\FF$
consider the cubic form
$$\6c_{t}:=(x_{1}^{(2)}+x_{2}^{(2)})y_{1}+x_{1}x_{2}y_{2}+
tx_{2}^{(2)}y_{3}$$ on $W_{1}\oplus W_{2}$, which is non-degenerate if
$t\ne0$.  A simple computation reveals that every $\6c_{t}$ is
contained in $\6A=\C\cap\6C_{q}$. The corresponding nil-polynomial
(depending on the choice of $\epsilon$) then is
$$\6p_{t}=\6q+\6c_{t}+\6d_{t}\Steil{with}\6d_{t}:=x_{1}^{(4)}
+x_{1}^{(2)}x_{2}^{(2)}+(1+\epsilon^{-1}t^{2}) x_{2}^{(4)}\,.\Leqno{GR}$$ 
For every $t\in\FF^{*}$ an invariant of $\6d_{t}$ is
$\pphi(t):=\qu{g_{2}(\6d_{t})^{3}}{g_{3}(\6d_{t})^{2}}=
\epsilon^{2}t^{-4}(4+\epsilon^{-1}t^{2})^{3}\in\FF$, where
$g_{2},g_{3}$ are the classical invariants of binary quartics, compare
\Lit{MUKA} \p27.  Since every fiber of $\pphi:\FF^{*}\to\FF$ contains
at most 6 elements we conclude

\Proposition{} For every field $\FF$ and every fixed
$\epsilon\in\FF^{*}$ the set of all equivalence classes given by all
nil-polynomials $\6p_{t}$, $t\in\FF^{*}$, has the same cardinality as
$\FF$ and, in particular, is infinite.  \Formend

{\medskip\noindent\bf Remarks} 1. In case $\FF=\QQ$ is the rational
field there are infinitely many choices of $\epsilon\in\QQ^{*}$
leading to pairwise non-equivalent quadratic forms $\6q$ in \Ruf{KU}.
For each such choice there is an infinite number of pairwise
non-equivalent nil-polynomials $\6p_{t}$ of degree 4 over $\QQ$.

\noindent 2. In case $\FF=\RR$ is the real field there are essentially
the two choices $\epsilon=\pm1$. In case $\epsilon=1$ the form $\6q$
has type $(5,2)$ and all nil-polynomials $\6p_{t}$ with $0<
t\le\sqrt8$ are pairwise non-equivalent.  In case $\epsilon=-1$ the
form $\6q$ has type $(4,3)$ and all $\6p_{t}$ with $t>0$
are pairwise non-equivalent.

\noindent 3. Nil-polynomials of degree $\ge5$ can be constructed just
as in the case of degrees 3 and 4 as before. As an example we briefly
touch the case of degree 5: Fix a vector space $W$ of finite dimension
over $\FF$ together with a non-degenerate quadratic form $\6q$ on
$W$. Assume furthermore that there is a direct sum decomposition
$W=W_{1}\oplus W_{2}\oplus W_{3}\oplus W_{4}$ into non-zero totally
isotropic subspaces such that $W_{1}\oplus W_{4}$ and $W_{2}\oplus
W_{3}$ are orthogonal. Then consider a non-degenerate cubic form $\6c$
on $W_{1}\oplus W_{2}\oplus W_{3}$ (trivially extended to $W$) that
can be written as a sum $\6c=\6c'+\6c''$ of cubic forms with the
following properties: $\6c'$ is a cubic form on $W_{1}\oplus W_{3}$
that is linear in the variables of $W_{3}$ while $\6c''$ is a cubic
form on $W_{1}\oplus W_{2}$ that is linear in the variables of
$W_{1}$. Denote by $x\cd y$ the commutative product on $W$ determined
by $\6q$ and $\6c$. Then $W_{2}=\langle W_{1}\cd W_{1}\rangle$,
$\,W_{3}=\langle W_{1}\cd W_{2}\rangle$ and $W_{4}=\langle W_{1}\cd
W_{3}+W_{2}\cd W_{2}\rangle$. If we assume $\6c\in\6C_{\6q}$, we get
with $W_{5}:=\FF$ as in the nil-index 4 case above that
$\5N:=\bigoplus_{k>0}W_{k}$ is a canonically graded admissible algebra
of nil-index 5. Conversely, every canonically graded admissible
algebra of nil-index 5 is obtained this way.

\KAP{counter}{Some counterexamples}  

In this section we give examples of admissible algebras {\sl without}
Property (AH). By Theorem \ruf{NB} such an algebra cannot be
graded. We also give examples of non-gradable algebras {\sl with}
Property (AH).

A test for non-gradability (especially with computer aid) involving
nil-polynomials is the following

\Proposition{UG} Let $\5N$ be an admissible algebra with
nil-polynomial $\6p\in\FF[x_{1},\dots,x_{n}]$ and nil-index
$\nu$. Then, if $\5N$ has a grading, there exists a matrix
$A=(a_{jk})\in\FF^{n\times n}$ such that \0 $\xi\6p=\6p$ for
$\xi:=\sum_{j,k=1}^{n}a_{jk}x_{k}\dd{x_{j}}$. \1 $A$ is diagonalizable
over $\FF$. \1 Every eigen-value of $A$ is a positive rational number,
and for every eigen-value $\epsilon$ also $1{-}\epsilon$ is an
eigen-value with the same multiplicity.\1 The eigenvalues form the
arithmetic progression ${1\over\nu},{2\over\nu},\dots,{\nu-1\over\nu}$
if $\5N$ has a canonical grading.

\Proof Let $\5N=\bigoplus\5N_{k}$ be a grading and put
$d:=\max\{k:\5N_{k}\ne0\}$. Then $\5N_{d}$ is the annihilator of
$\5N$. Let $W:=\bigoplus_{k<d}\5N_{k}$ and choose a linear isomorphism
$\phi:\FF^{n}\to W$ such that $\phi(e_{j})\in\5N_{k}$ for all $j$ and
suitable $k=k(j)$, where $e_{1},\dots,e_{n}$ is the standard basis of
$\FF^{n}$. Further let $\omega$ be a pointing on $\5N$ with kernel $W$
and put $\tilde\6p:=\omega\circ\exp_{2}\circ\,\phi$. Then (i) - (iii)
hold for $\tilde\6p$ and the diagonal matrix $\tilde A$ with diagonal
entries $\tilde a_{jj}=k(j)/d$ for all $j$. In case the grading is
canonical, $d=\nu$ and $\5N_{k}\ne0$ for $0<k<\nu$ holds, implying
(iv).

By Proposition \ruf{UQ} the nil-polynomials $\6p$, $\tilde\6p$ are
linearly equivalent. As a consequence of Proposition \ruf{HD} there
exists $C\in\GL(n,\FF)$ and $c\in\FF^{*}$ with
$c\,\tilde\6p(x)=\6p(Cx)$ for all $x\in\FF^{n}$. But then $\6p$
satisfies (i) - (iv) with respect to $A=C\tilde AC^{-1}$.\qed

\smallskip The proof of \ruf{UG} uses the non-trivial fact that every
gradable $\5N\,$ has property (AH) and hence that any two associated
nil-polynomials are linearly equivalent, compare \ruf{NB} and
\ruf{NE}.  Property (i) means that
$\6p=\lambda_{1}\!\raise1.5pt\hbox{$\,\partial
\6p$}/\raise-2.5pt\hbox{$\!\partial
x_{1}\,$}+\dots+\lambda_{n}\!\raise1.5pt\hbox{$\,\partial
\6p$}/\raise-2.5pt\hbox{$\!\partial x_{n}\,$}$ for suitable linear
forms $\lambda_{k}$ on $\FF^{n}$, and hence that $\6p$ is in its
Jacobi ideal (the ideal in $\FF[x_{1},\dots,x_{n}]$ generated by all
first partial derivatives of $\6p$).  Part of \ruf{UG} can also be
reformulated in terms of quasi-homogeneous polynomials. By definition,
$f\in\FF[x_{1},\dots,x_{n}]$ is {\sl quasi-homogeneous} if there exist
positive integers $m,m_{1},\dots,m_{n}$ with
$f(t^{m_{1}}x_{1},\dots,t^{m_{n}}x_{n})\equiv
t^{m}f(x_{1},\dots,x_{n})$ for all $t\in\FF$.

\noindent{\bf Remark.} Let $\5N$ be an admissible algebra having a
grading. Then there exists a quasi-homoge\-neous nil-polynomial
associated with $\5N$.\Formend

\medskip Motivated by the main result of \Lit{SAIT}, see also
\Lit{XUYA}, we organize our search for non-gradable admissible
algebras as follows: For fixed indeterminates $T_{1},\dots,T_{m}$
denote by $\7m$ the maximal ideal in the localization $R$ of
$\FF[T_{1}\dots,T_{m}]$ at the origin, that is, the ideal of all
quotients $P/Q$ with $P,Q\in\FF[T_{1}\dots,T_{m}]$ satisfying $Q(0)\ne
P(0)=0$. For every $F\in\7m^{2}$ let $J(F)$ in $R$ be the {\sl
Jacobi ideal} of $F$ (the ideal generated by all first order
derivatives of $F$). Then if $J(F)\supset\7m^{k}$ for some $k$,
the Milnor algebra $R/J(F)$ has finite dimension and its maximal ideal
$$\5M(F):=\7m/J(F)\Leqno{PQ}$$ is nilpotent. For our search we are
looking for examples with $F\notin J(F)$.

Let us stress that we always use the notation $\5M(F)$ in the
following way: The indeterminates $T_{1},\dots,T_{m}$ giving the
localization $R$ above are precisely those occurring in $F$. As an
example, for the $\5M(F)$ occurring two lines below we understand
$m=2$ and $\{T_{1},T_{2}\}=\{X,Y\}$.

\bigskip \medskip \centerline{\bf Non-gradable algebras with Property
(AH)}

Let $\5N:=\5M(X^{5}+X^{2}Y^{2}+Y^4)$. Then $\5N$ has basis
$$x,x^2,x^3,x^4,y,xy,y^2,y^3,x^5,\Leqno{UB}$$ where $x,y$ are the
residue classes of $X,Y$. We abbreviate this basis with
$e_{1},\dots,e_{9}$. The annihilator $\5N_{[1]}$ is spanned by
$e_{9}=x^{5}$ and the residue class of $F$ is
$-e_{9}/4\in\5N_{[1]}\,$.  A nil-polynomial
$\6p\in\FF[x_{1},\dots,x_{8}]$ obtained from the basis
is{\parindent10pt\footnote{$^{1}$}{\sevenrm Computed with {\eightt
Singular}, freely available at {\eightt
http://www.singular.uni-kl.de/}}}
$$\eqalign{
&\textstyle{1\over120}x_{1}^5+\big({1\over6}x_{1}^3x_{2}^{}+
{5\over96}x_{5}^4-{5\over8}x_{1}^2x_{5}^2\big)\cr
&\textstyle\qquad+\big({1\over2}x_{1}^2x_{3}^{}+{1\over2}x_{1}^{}x_{2}^2+
{5\over8}x_{5}^2x_{7}^{}-{5\over4}x_{1}^2x_{7}^{}-
{5\over2}x_{1}^{}x_{5}^{}x_{6}^{}-{5\over4}x_{2}^{}x_{5}^2\big)\cr
&\textstyle\qquad+\big(x_{1}^{}x_{4}^{}+x_{2}^{}x_{3}^{}+
{5\over4}x_{5}^{}x_{8}^{}+{5\over8}x_{7}^2-
{5\over2}x_{2}^{}x_{7}^{}-{5\over4}x_{6}^2\big)\;.\cr }\Leqno{UV}$$
\vskip-8pt\noindent In particular, $\5N$ is an admissible algebra of
dimension 9 and nil-index 5. Also, $\5N$ does not have a gradation but
has Property (AH){\parindent10pt\footnote{$^{2}$}{\sevenrm
Computed with {\eightt Maple.}}}.

Further examples of this type (but of higher dimension) are obtained
by varying $F$. For instance $\5N:=\5M(F)$ with
$F=X^{4}+X^{2}Y^{3}+Y^{5}$ is an admissible algebra of dimension 11
and nil-index 5 without a grading but with Property (AH). As algebra
$\5N$ is isomorphic to the (unique) maximal ideal of
$\FF[X,Y]/I$, where $I:=\big(\raise1.5pt\hbox{$\,\partial
F$}/\raise-2.5pt\hbox{$\!\partial X\,$},\raise1.5pt\hbox{$\,\partial
F$}/\raise-2.5pt\hbox{$\!\partial Y\,$},\,X^{3}Y\big)$. The latter
algebra already appears as an example of a non-gradable algebra in
Remark 3.3 in \Lit{COKR} (note that $X^{5}$ occurring there is already
contained in $I$ and hence is superfluous).

If we go to algebras with embedding dimension 3 we can get an
example with nil-index 4 and dimension 8. Add an additional
indeterminate $Z$ and consider $\5M(X^{4}+ XY^{2}+Y^{3}+XZ^{2})$. Then
a basis is given by $x,y,z,x^2,x^3,yz,z^2,x^4$ and
$$\eqalign{\6p=
&\textstyle{1\over24}x_{1}^4+({1\over2}x_{1}^2x_{4}^{}+
3x_{1}^2x_{2}^{}-2x_{1}^{}x_{2}^2+{4\over9}x_{2}^3
-{4\over3}x_{2}^{}x_{3}^2)\cr
&\textstyle\qquad+(x_{1}^{}x_{5}^{}+{1\over2}x_{4}^2+6x_{2}^{}x_{4}^{}
-{8\over3}x_{2}^{}x_{7}^{}-{8\over3}x_{3}^{}x_{6}^{})\cr }$$ 
is the
nil-polynomial derived from it. The Hilbert function is
$\{1,3,3,1,1\}$. Also for this algebra there is no matrix
$A\in\FF^{7\times7}$ satisfying (i) in Proposition \ruf{UG}. On the
other hand, for suitable coefficients
$\lambda_{k}\in\FF[x_{1},\dots,x_{7}]$ there exists$^{1}$ a
representation $\6p=\lambda_{1}\!\raise1.5pt\hbox{$\,\partial
\6p$}/\raise-2.5pt\hbox{$\!\partial
x_{1}\,$}+\dots+\lambda_{7}\!\raise1.5pt\hbox{$\,\partial
\6p$}/\raise-2.5pt\hbox{$\!\partial x_{7}\,$}$, that is, $\6p$ lies in
its Jacobi ideal.

\vfill\eject
\bigskip\centerline{\bf Failing Property (AH)}

In this subsection we restrict to the special case
$\FF\in\{\RR,\CC\}$. This allows us to use analytic arguments.

It can be seen$^{1}$ that $\5M(X^{6}+X^{2}Y^{3}+Y^{5})$ is an
admissible algebra of dimension 17 and nil-index 7.  In the same way
$\5M(X^{7}+X^{2}Y^{3}+X^{3}Y^{2}+Y^{4})$ is an admissible algebra of
dimension 15 and nil-index 8. Also
$\5M(X^{5}+X^{2}Y^{2}+Y^{4}+XZ^{2})$ for fixed indeterminates $X,Y,Z$
is an admissible algebra of dimension 15 with nil-index 6.  All three
algebras do not have Property (AH). The automorphism groups have
dimension $25,20,23$ respectively.  Instead of giving further details
for these examples we do this for another one (of even lower nil-index
but higher dimension).

For indeterminates $X,Y,Z,U$ consider
$\5N:=\5M(X^{3}+X^{2}Y^{2}+Y^{4}+XZ^{2}+ZU^{2})$. Then $\5N$ is an
admissible algebra of dimension 23 and nil-index 5, in fact, $\5N$ has
symmetric Hilbert function $\{1,4,7,7,4,1\}$. A basis for $\5N$
is$^{1}$ \medskip \centerline{\klein
$x,y,z,u,xy,y^2,xy^2,yz,y^2\!z,z^2,yz^2,y^2\!z^2,xu,yu,
xyu,y^2\!u,xy^2\!u,u^2,xu^2,yu^2,xyu^2,y^2\!u^2,xy^2\!u^2 $}
\medskip\noindent with the last vector spanning $\Ann(\5N)$. Also, the
nil-polynomial given by the above basis is$^{1}$

\vskip-13pt
$${\klein\eqalign{\6p=
&\textstyle~\big({1\over4}x_{1}^{}x_{4}^2-{1\over8}x_{1}^2x_{3}^{}+{1\over96}x_{2}^2x_{3}^{}+{1\over8}x_{3}^3\big)x_{2}^2\cr
&\textstyle+\big({1\over18}x_{1}^3x_{3}^{}-{1\over4}x_{1}^2x_{2}^{}x_{8}^{}-{1\over4}x_{1}^2x_{6}^{}x_{3}^{}-{1\over6}x_{1}^2x_{4}^2+{1\over2}x_{1}^{}x_{2}^2x_{18}^{}-{1\over2}x_{1}^{}x_{2}^{}x_{5}^{}x_{3}^{}+{1\over4}x_{6}^{}x_{3}^3+x_{1}^{}x_{2}^{}x_{4}^{}x_{14}^{}\cr
&\textstyle\quad+{1\over2}x_{1}^{}x_{6}^{}x_{4}^2+{1\over24}x_{2}^3x_{8}^{}+{1\over8}x_{2}^2x_{6}^{}x_{3}^{}+{3\over4}x_{2}^2x_{3}^{}x_{10}^{}+{1\over2}x_{2}^2x_{4}^{}x_{13}^{}+{1\over2}x_{2}^{}x_{5}^{}x_{4}^2+{3\over4}x_{2}^{}x_{3}^2x_{8}^{}\big)\cr
&\textstyle+\big(x_{1}^{}x_{2}^{}x_{20}^{}-{1\over4}x_{1}^2x_{9}^{}-{1\over3}x_{1}^2x_{18}^{}-{1\over2}x_{1}^{}x_{5}^{}x_{8}^{}+x_{1}^{}x_{6}^{}x_{18}^{}-{1\over2}x_{1}^{}x_{7}^{}x_{3}^{}-{2\over3}x_{1}^{}x_{4}^{}x_{13}^{}+x_{1}^{}x_{4}^{}x_{16}^{}\cr
&\textstyle\quad+x_{2}^{}x_{13}^{}x_{14}^{}+{1\over8}x_{2}^2x_{9}^{}+{1\over2}x_{2}^2x_{19}^{}+x_{2}^{}x_{5}^{}x_{18}^{}+{1\over4}x_{2}^{}x_{6}^{}x_{8}^{}+{3\over2}x_{2}^{}x_{3}^{}x_{11}^{}+{3\over2}x_{2}^{}x_{8}^{}x_{10}^{}+x_{2}^{}x_{4}^{}x_{15}^{}\cr
&\textstyle\quad+{1\over2}x_{1}^{}x_{14}^2-{1\over4}x_{5}^2x_{3}^{}+x_{5}^{}x_{4}^{}x_{14}^{}+{1\over8}x_{6}^2x_{3}^{}+{3\over2}x_{6}^{}x_{3}^{}x_{10}^{}+x_{6}^{}x_{4}^{}x_{13}^{}+{1\over2}x_{7}^{}x_{4}^2+{3\over4}x_{3}^2x_{9}^{}+{3\over4}x_{3}^{}x_{8}^2\big)\cr
&\textstyle+\big(x_{1}^{}x_{22}^{}-{2\over3}x_{1}^{}x_{19}^{}+x_{2}^{}x_{21}^{}+x_{5}^{}x_{20}^{}+{1\over4}x_{6}^{}x_{9}^{}+x_{6}^{}x_{19}^{}+x_{7}^{}x_{18}^{}+{3\over2}x_{3}^{}x_{12}^{}\cr
&\textstyle\quad+{3\over2}x_{8}^{}x_{11}^{}+{3\over2}x_{9}^{}x_{10}^{}+x_{4}^{}x_{17}^{}-{1\over3}x_{13}^2+x_{13}^{}x_{16}^{}+x_{14}^{}x_{15}^{}\big)\;.\cr
}}\Leqno{SR}$$

\vskip-8pt\noindent We claim that the graph $S$ of $\6p$ in $\FF^{23}$
cannot be affinely homogeneous. Define $\6f$ on $\FF^{23}$ by
$\6f(x_{1},\dots,x_{23}):=\6p(x_{1},\dots,x_{22})-x_{23}$. Since
$\Aff(S)$ is a Lie group over $\FF$, affine homogeneity of $S$ would
imply for every fixed $1\le\ell\le22$ the existence of an affine
vector field
$$\xi=\dd{x_{\ell}}+\sum_{j,k=1}^{23}a^{}_{jk}x^{}_{k}\!\dd{x_{j}}$$
with coefficients $a_{jk}\in\FF$ such that $\xi\6f=\rho\6f$ for some
$\rho\in\FF$. Since the polynomial $\6f$ has rational coefficients the
vector field $\xi$ can be chosen in such a way that all $a_{jk}$ and
$\rho$ are rational.  But $\xi\6f=\rho\6f$ is equivalent to a system
of linear equations in $a_{jk},\rho$. It can be seen$^{2}$ that this
system has a rational solution only in case $\ell\ne3$. Therefore the
orbit of $0\in S$ under the group $\Aff(S)$ has dimension 21. In
particular, $S$ is not even locally affinely homogeneous at the
origin.  Furthermore$^{2}$, the group $\Aff(S)\cong\Aut(\5N)$ is a Lie
group of dimension $42$ over $\FF$.

\bigskip\noindent{\bf Remark.} Since computer software may contain
errors we carried out all computations in a highly redundant manner:
Different machines were used, the same computations were repeated with
different routines, bases of the admissible algebras $\5N$ under
consideration were changed resulting in totally different (but
affinely equivalent) nil-polynomials and finally, the linear system
$\xi\6f=\rho\6f$ with $\xi$ from above was replaced by a totally
different one (namely Lemma \ruf{JP} for the special case $r=22$ and
$n=23$). Notice that $\xi\6f{-}\rho\6f\in\FF[x_{1},\dots,x_{22}]$ is a
polynomial of degree 5 having $a_{jk}$ and $\rho$ as parameters. An
alternative condition for $\xi$ being tangent to $S\subset\5N$ clearly
is that $\xi\6f$ vanishes at every point of $S$. This leads to a
linear system in terms of a polynomial of degree bigger than 5. This
works also for arbitrary codimension:

\medskip For fixed integers $r,s\ge1$ and $n:=r+s$ consider $\FF^{r}$,
$\FF^{s}$ with coordinates $(x_{1},\dots,x_{r})$,
$(x_{r+1},\dots,x_{n})$ respectively. Let $\6p:\FF^{r}\to\FF^{s}$ be a
polynomial map. Then $\6p=(\6p_{r+1},\dots,\6p_{n})$ for scalar valued
polynomials $\6p_{k}$ on $\FF^{r}$. Denote by $S$ the graph of $\6p$
in $\FF^{n}$, which is a smooth variety of codimension $s$. The proof
of the following statement is an easy exercise in differentiation and
will be omitted.

\Lemma{JP} Let $\FF$ be either $\RR$ or $\CC$. Then $S\subset\FF^{n}$
is locally affinely homogeneous at the origin if and only if for every
$1\le\ell\le r$ the following linear system in the unknowns
$\,a_{jk}$, $\,1\le j,k\le n$, has a solution over $\FF$
$$\6P_{m}=\sum_{j=1}^{r}(\6P_{j}+\delta_{j\ell}\,)\dD{\6p_{m}}{x_{j}}
\Steil{for}m=r+1,\dots,n\,,\;\hbox{where}\leqno{(\dagger)}$$
$$\6P_{j}:=\sum_{k=1}^{r}a_{jk}x_{k}+\sum_{k=r+1}^{n}a_{jk}\6p^{}_{k}
\Steil{for}1\le j\le n$$ and $\delta$ is the Kronecker delta.
\Formend

\medskip The equations in Lemma \ruf{JP} can also be used to compute
the Lie algebra $\aff(S)$ of $\Aff(S)$ numerically. Indeed, add to the
$a_{jk}$ further unknowns $c_{1},\dots,c_{r}$ and replace in
$(\dagger)$ the term $\delta_{j\ell}$ by $c_{j}$. Then the solution
space for this altered linear system is canonically isomorphic to
$\aff(S)$.

\bigskip\bigskip

\vskip7mm {\gross\noindent References} \medskip {\klein\openup-.5pt
\parindent 15pt\advance\parskip-1pt

\Ref{COKR}Corti\~nas, G., Krongold, F.: Artinian algebras and differential forms. Comm. Algebra {\bf 27}, 1711 - 1716 (1999).
\Ref{EAWO}Eastwood, M.G.: Moduli of isolated hypersurface singularities. Asian J. Math. {\bf 8}, 305-314 (2004). 
\Ref{ELRO}Elias, J., Rossi, M.E.: Isomorphism classes of short Gorenstein local rings via Macaulay's inverse system. TAMS to appear.
\Ref{FKAP}Fels, G., Kaup, W.: Local tube realizations of CR-manifolds
and maximal abelian subalgebras. Ann. Sc. Norm. Sup. Pisa. Cl. Sci. {\bf X}, 99-128 (2011).
\Ref{FEKA}Fels, G., Kaup W.: Classification of commutative algebras and tube realizations of hyperquadrics.\nline arXiv:0906.5549v2.
\Ref{FIKK}Fels, G., Isaev, A., Kaup W., Kruzhilin N.: Singularities and Polynomial Realizations of Affine Quadrics. J. Geom. Analysis {\bf 21} (2011).
\Ref{ISAA}Isaev, A.V.: On the number of affine equivalence classes of spherical tube hypersurfaces. Math. Ann. {\bf 349}, 59-72 (2011).
\Ref{ISAE}Isaev, A.V.: On the Affine Homogeneity of Algebraic Hypersurfaces Arising from Gorenstein Algebras.\nline a ~http://arxiv.org/pdf/1101.0452v1.\nline b ~http://arxiv.org/pdf/1101.0452v2.
\Ref{MUKA}Mukai, S.: An introduction to invariants and moduli. Cambridge Univ. Press 2003.
\Ref{PERE}Perepechko, A.: On solvability of the automorphism group of a finite-dimensional algebra. arXiv:1012:0237
\Ref{SAIT}Saito, K.: Quasihomogene isolierte Singularit\"aten von Hyperpl\"achen. Invent. Math. {\bf14}, 123-142 (1971).
\Ref{XUYA}Xu, Y.J., Yau, S.S.T.: Micro-local characterization of quasi-homogeneous singularlities. Amer. J., Math. {\bf118}, 389-399 (1996).

\bigskip
}

\smallskip\openup-4pt\parindent0pt
\hbox to 3cm{\hrulefill}

\smallskip

Universit\"at T\"ubingen, Mathematisches Institut

Auf der Morgenstelle 10

D-72076 T\"ubingen, Germany 

\medskip
e-mail: {\tt gfels@uni-tuebingen.de}

e-mail: {\tt kaup@uni-tuebingen.de\par\vskip8pt}

\closeout\aux\bye